\documentclass[a4,12pt]{article}
\usepackage[centertags]{amsmath}
\usepackage{amsmath}
\usepackage{amssymb}
\usepackage{amsthm}
\usepackage{color}
\usepackage{dsfont}
\usepackage{latexsym}
\usepackage[T1]{fontenc}
\usepackage[latin1]{inputenc}
\usepackage{fancyhdr}
\usepackage{indentfirst}
\usepackage{times}
\usepackage{enumerate}
\usepackage{cite}
\usepackage{ragged2e}
\usepackage[colorlinks=true,linkcolor=blue,citecolor=red,urlcolor=black]{hyperref}
\newtheorem{definition}{Definition}[]

\newtheorem{theorem}{Theorem}[]
\newtheorem{proposition}{Proposition}[]
\newtheorem{lemma}{Lemma}[]
\newtheorem{remark}{Remark}[]
\newtheorem{example}{Example}[]
\newtheorem{corollary}{Corollary}[]
\def \ind{1\!\!1}
\title{Deterministic Differential Games in Infinite Horizon Involving Continuous and Impulse Controls}
\author{Brahim El Asri\thanks{Ibn Zohr University, Laboratoire LISAD, Equipe Aide à la Decision, ENSA, B.P. 1136, Agadir, Morocco, e-mail: b.elasri@uiz.ac.ma.} \, and\, Hafid Lalioui\thanks{Ibn Zohr University, Laboratoire LISAD, Equipe Aide à la Decision, ENSA, B.P. 1136, Agadir, Morocco, e-mail: hafid.lalioui@edu.uiz.ac.ma. Financially supported by CNRST, Rabat, Morocco (Grant 17 UIZ 2019).}}

\begin{document}

\date{}\maketitle

\begin{abstract}

We study a new class of two-player, zero-sum, deterministic, differential games where each player uses both continuous and impulse controls in infinite horizon with discounted payoff. We assume that the form and cost of impulses depend on nonlinear functions and on the state of the system, respectively. We use the Bellman's dynamic programming principle (DPP) and viscosity solutions approach to show, for this class of games, existence and uniqueness of a solution for the associated Hamilton-Jacobi-Bellman-Isaacs (HJBI) partial differential equations (PDEs). We then, under Isaacs' condition, deduce that the lower and the upper value functions coincide.

\end{abstract}

\noindent
\textbf{Keywords:} Deterministic differential game, Infinite horizon, Continuous and impulse control, Dynamic programming principle, Viscosity solution, Isaacs' condition.\\\\
\noindent
\textbf{AMS Classifications (2010):} 49K35, 49L25, 49N70, 90C39, 93C20.

\section{Introduction}

We consider a new class of deterministic infinite horizon two-player zero-sum differential games where both continuous and impulse controls are involved. Deterministic differential games with continuous controls alone, started with the work of Pontryagin \& al \cite{BGMP62} and Isaacs \cite{I65}, have been studied in Bardi and Capuzzo-Dolcetta \cite{BC97} and in Evans and Souganidis \cite{ES84} with non-anticipative strategy in the viscosity solutions framework. Zero-sum Differential games with impulse controls were studied in Yong \cite{Y94} for deterministic case with only one impulse control, in Barles \cite{B85} and El Farouq \& al \cite{ElBB10} the authors considered continuous and one impulse controls problem, and recently in Cosso \cite{C13}, Azimzadeh \cite{A17} and El Asri and Mazid \cite{ElM18} zero-sum stochastic games with two impulse controls were studied (see also Zhang \cite{Z11}). In the previous literature of differential games one can find those of mixed type, Dharmatti \& al \cite{DMM09}, where the state is controlled by a combination of both continuous and discrete actions, and those with continuous, switching, hybrid and impulse controls (Dharmatti and Shaiju \cite{DS05,DS07}, Dharmatti and Ramaswamy \cite{DR05,DR06} and Barles \& al \cite{BDR10}). Let us just mention that in Yong \cite{Y90,Y'90} differential games where both players use switching controls are studied. In \cite{Y94}, where zero-sum differential games involving impulse controls are considered, one player is using continuous control whereas the opponent uses impulse control. We also mention that \cite{DS05} extends the work in \cite{Y94} to a two-person zero-sum differential game involving continuous, switching and impulse controls, where the existence of value and its characterization as the unique viscosity solution of the associated system of quasi-variational inequalities (QVIs) have been proved.
\par In this paper, we study a deterministic infinite horizon two-player zero-sum differential game where the two players adopt both continuous and impulse controls, where the form and cost of impulses depend on nonlinear functions and on the state of the system, respectively. The state $y_{x}(t)$ of the continuous and impulse controls game considered evolves according to the following dynamical system:
\begin{equation}\label{system}\tag{S}
\left\{
\begin{aligned}
y_x(0^-)&=x\in\mathbb{R}^n\;\text{(initial state)};\\
\dot y_x(t)&=b\bigl(y_x(t);\theta_1(t),\theta_2(t)\bigr),\;t\in\mathbb{R}^+,\;t\neq\tau_m\;\text{and}\;t\neq\rho_k;\\
y_x(\tau_m^+)&=y_x(\tau_m^-)+g_\xi\bigl(y_x(\tau_m^-),\xi_m\bigr)\prod_{k\geq 0}\ind_{\{\tau_m\neq\rho_k\}},\;\tau_m\geq 0,\;\xi_m\neq 0;\\
y_x(\rho_k^+)&=y_x(\rho_k^-)+g_\eta\bigl(y_x(\rho_k^-),\eta_k\bigr),\;\rho_k\geq 0,\;\eta_k\neq 0.
\end{aligned}
\right.
\end{equation}
Here $b$ is a function from $\mathbb{R}^n\times\mathbb{R}^l\times\mathbb{R}^l$ into $\mathbb{R}^n$, $\theta_1:t\rightarrow\theta_1(t)$ $\bigl(\text{resp.}\;\theta_2:t\rightarrow\theta_2(t)\bigr)$ being the continuous control, is any function from $\Theta_1$ $(\text{resp.}\;\Theta_2)$ the space of measurable functions from $\mathbb{R}^+$ to a compact metric space $A$ (resp. $B$) subset of $\mathbb{R}^l$. The two functions $g_\xi$ and $g_\eta$ are from $\mathbb{R}^n\times\mathbb{R}^p$ into $\mathbb{R}^n$ and $\mathbb{R}^n\times\mathbb{R}^q$ into $\mathbb{R}^n$, respectively. The sequences $\bigl(\{\tau_m\},\{\xi_m\}\bigr)_{m\in\mathbb{N}}$ and $\bigl(\{\rho_k\},\{\eta_k\}\bigr)_{k\in\mathbb{N}}$ represent the two impulse controls, where $\{\tau_m\}_{m\in\mathbb{N}}$ and $\{\rho_k\}_{k\in\mathbb{N}}$ are two non-decreasing sequences of non-negative real numbers which satisfy:
$$\tau_m\rightarrow+\infty\;\text{when}\;m\rightarrow+\infty,\;\text{and}\;\rho_k\rightarrow+\infty\;\text{when}\;k\rightarrow+\infty,$$
and $\{\xi_m\}_{m\in\mathbb{N}}$ and $\{\eta_k\}_{k\in\mathbb{N}}$ are two sequences of elements of convex cones $U\subset\mathbb{R}^p$ and $V\subset\mathbb{R}^q$, respectively. We mention that the state value $y_x(t)$ of the controlled system at time $t$ is driven by the two continuous controls $\theta_1(.)\in\Theta_1$ and $\theta_2(.)\in\Theta_2$ for $player-\xi$ and $player-\eta$, respectively. In addition, both players are allowed to use impulse controls, $u$ for $player-\xi$ and $v$ for $player-\eta$. The impulse controls $u$ and $v$ are defined by the double sequences $u:=(\tau_m,\xi_m)_{m\in\mathbb{N}}$ and $v:=(\rho_k,\eta_k)_{k\in\mathbb{N}}$, respectively, where the actions $\xi_m$ and $\eta_k$ belong to the spaces of impulse control actions $U$ and $V$, respectively. The infinite product $\prod_{k\geq 0}\ind_{\{\tau_m\neq\rho_k\}}$ signifies that when the two players act together on the system at the same time using the impulse controls, we take into account only the action of $player-\eta$. The gain (resp. cost) functional (\ref{payoff}) for $player-\xi$ (resp. $player-\eta$) is defined as follows:
\begin{equation}\label{payoff}\tag{$J$}
\begin{aligned}
J\bigl(x;\theta_1(.),u,\theta_2(.),v\bigr):=&\int_0^\infty f\bigl(y_x(t);\theta_1(t),\theta_2(t)\bigr)\exp(-\lambda t)dt\\
&-\sum_{m\geq 0}c\bigl(y_x(\tau_m^-),\xi_m\bigr)\exp(-\lambda\tau_m)\prod_{k\geq 0}\ind_{\{\tau_m\neq\rho_k\}}\\
&+\sum_{k\geq 0}\chi\bigl(y_x(\rho_k^-),\eta_k\bigr)\exp(-\lambda\rho_k),
\end{aligned}
\end{equation}
where $y_x(t):=y_x\bigl(t;\theta_1(.),u,\theta_2(.),v\bigr)$, and $c$ and $\chi$ are two positive impulse cost functions for $player-\xi$ and $player-\eta$, respectively, which represent the cost of impulse actions for both players. We note that the cost of a player is the gain for the other (zero-sum), meaning that when a player performs an impulse action he/she has to pay a positive cost, resulting in a gain for the other player. The function $f$ from $\mathbb{R}^n\times\mathbb{R}^l\times\mathbb{R}^l$ into $\mathbb{R}$ represents the running gain (resp. cost) for maximizing $player-\xi$ $(\text{resp. minimizing}\;player-\eta)$ and the positive real constant $\lambda$ being the discount factor.
\par To the best of our knowledge, the large literature on differential games involving continuous and/or impulse controls does not provide any theoretical and/or computational means to study the class of games we have considered in system (\ref{system}), related to functional (\ref{payoff}). In the framework of differential games the biggest difficulty lies in showing the characterization of the value function in the viscosity solutions sense. It is challenging to show the comparison theorem which is necessary to get uniqueness for the characterization. We mention that the following assumptions, usually used to deal with impulse control problems (see e.g. \cite{DS05,Y94,Z11}), for $x\in\mathbb{R}^n$, $\xi_1,\xi_2\in U$ and $\eta_1,\eta_2\in V$, $c(x,\xi_1+\xi_2)<c(x,\xi_1)+c(x,\xi_2)$ and/or $\chi(x,\eta_1+\eta_2)<\chi(x,\eta_1)+\chi(x,\eta_2),$ were dropped. Relative to a large part of the existing work, another contribution of the paper is the study of this class of differential games with form and cost of impulses depending, respectively, on the nonlinear functions $g_\xi,g_\eta$ and system's state $y_x(.)$. Works in references \cite{BDR10,DR05}, where a hybrid control system in which both discrete and continuous controls are involved was studied in the viscosity solutions framework, are our closest related papers (see also \cite{B85,ElBB10,DS05,DR06,DS07,DMM09,Y89,Y90,Y'90,Y94,Z11}).
Regarding recent related papers, Bettiol \& al \cite{BQV19} consider a two-player zero-sum differential game with a cost of Bolza type, subject to a state constraint, in El Farouq \cite{El17} the author has proved that the value function of a deterministic infinite horizon, mixed continuous and impulse control problem, is the unique viscosity solution of a related first order Hamilton-Jacobi quasi-variational inequality (QVI), and in A\"id \& al \cite{ABCCT20}, Basei \& al \cite{BCG21} and Sadana \& al \cite{SRZ21,SRB21} some nonzero-sum game problems with impulse controls were studied. In our previous work El Asri \& al \cite{ElLM21}, we have considered a zero-sum deterministic game problem where only impulse controls are involved in infinite-time horizon, where a new Hamilton-Jacobi-Bellman-Isaacs (HJBI) QVI was given to prove, under a proportional property assumption on the maximizer cost function, that the value functions coincide and turn out to be the unique viscosity solutions to the given HJBI QVI.
\par The terminology of a QVI, introduced to deal with impulse control problems in Bensoussan and Lions \cite{BL84}, the definition of lower and upper value functions for differential games, defined in Elliott and Kalton \cite{EK72,EK74} and in \cite{ES84}, and the viscosity solutions approach in Crandall and Lions \cite{CL83} and Crandall et al \cite{CEL84, CIL92}, lead to characterize the value of a game as the unique viscosity solution of its related QVI. The relationship between the two-player, zero-sum, deterministic, differential games and the theory of viscosity solutions was first shown in \cite{ES84}, Barron et al \cite{BEJ84} and Souganidis \cite{S85,S85'}. Our aim in this work lies in the fact that we investigate, via the theory of viscosity solutions, the class of two-player, zero-sum, deterministic, continuous and impulse controls games given by the system (\ref{system}), in infinite horizon. Indeed, we describe the problem by the following associated double-obstacle lower and upper HJBI equations (\ref{LowerQVI}) and (\ref{UpperQVI}), respectively, where the Hamiltonians $H^-$ and $H^+$ involve only the first order partial derivatives:
\begin{equation}\label{LowerQVI}\tag{L}
\min\biggl\{\max\Bigl[\lambda v(x)+H^-\bigl(x,Dv(x)\bigr),v(x)-\mathcal{H}_{inf}^\chi v(x)\Bigr],v(x)-\mathcal{H}_{sup}^c v(x)\biggr\}=0,
\end{equation}
and
\begin{equation}\label{UpperQVI}\tag{U}
\max\biggl\{\min\Bigl[\lambda v(x)+H^+\bigl(x,Dv(x)\bigr),v(x)-\mathcal{H}_{sup}^c v(x)\Bigr],v(x)-\mathcal{H}_{inf}^\chi v(x)\biggr\}=0,
\end{equation}
where $Dv(.)$ denotes the gradient of the function $v:\mathbb{R}^n\rightarrow\mathbb{R}$, $v$ represents the value function of the game problem, the lower Hamiltonian $H^-$ is given by:
$$H^-\bigl(x,Dv(x)\bigr):=\inf_{\theta_1\in A}\sup_{\theta_2\in B}\bigl(-Dv(x).b(x;\theta_1,\theta_2)-f(x;\theta_1,\theta_2)\bigr),$$
and the upper Hamiltonian $H^+$ is defined as follows:
$$H^+\bigl(x,Dv(x)\bigr):=\sup_{\theta_2\in B}\inf_{\theta_1\in A}\bigl(-Dv(x).b(x;\theta_1,\theta_2)-f(x;\theta_1,\theta_2)\bigr).$$
The first (resp. second) obstacle is defined through the use of the minimum (resp. maximum) cost operator $\mathcal{H}_{inf}^\chi$ (resp. $\mathcal{H}_{sup}^c$), where
$$\mathcal{H}_{inf}^\chi v(x):=\inf_{\eta\in V}\Bigl[v\bigl(x+g_\eta(x,\eta)\bigr)+\chi(x,\eta)\Bigr]\;\Bigl(\text{resp}.\;\mathcal{H}_{sup}^c v(x):=\sup_{\xi\in U}\Bigl[v\bigl(x+g_\xi(x,\xi)\bigr)-c(x,\xi)\Bigr]\Bigr).$$
\par Our main results concern the existence and uniqueness of solution in viscosity sense for both HJBI QVIs (\ref{LowerQVI}) and (\ref{UpperQVI}). More specifically, under classical assumptions of the impulse control problems using the Bellman's dynamic programming principle (DPP) for the class of differential games studied, we prove that the lower and the upper value functions are viscosity solutions of the lower HJBI QVI (\ref{LowerQVI}) and the upper HJBI QVI (\ref{UpperQVI}), respectively. Moreover, by reformulating these QVIs, we establish a Comparison Theorem which stands for the major contribution of the paper and shows the uniqueness results in the viscosity solutions sense for these QVIs. Next, we state the Isaacs' condition $H^-=H^+$ for the differential game to have a value.
\par The outline of the paper is the following: in Section \ref{Sect.2}, we present the infinite horizon, zero-sum, deterministic, continuous and impulse controls game studied and we give its related definitions and assumptions. Further, we give regularity results for the associated lower and upper value functions of the game, we show first that both satisfy the DPP property, then we prove that they are bounded and uniformly continuous (BUC) in $\mathbb{R}^n$. Section \ref{Sect.3} is devoted to the viscosity characterization of the corresponding lower and upper HJBI QVIs, it is shown that the lower and the upper value functions are viscosity solutions to the lower HJBI QVI (\ref{LowerQVI}) and the upper HJBI QVI (\ref{UpperQVI}), respectively. In Section \ref{Sect.4}, we establish the Comparison Theorem which gives the uniqueness result for both lower and upper HJBI QVIs. Hence, under Isaacs' condition, we deduce that the game admits a value.

\section{Statement of the Continuous and Impulse Controls Differential Game and Classic Results}\label{Sect.2}

\subsection{Zero-Sum Deterministic Continuous and Impulse Controls Game}\label{Sect.2.1}

We consider the two-player, zero-sum, deterministic, differential game in $\mathbb{R}^n$ described, in the introduction show, by the system (\ref{system}) and the functional (\ref{payoff}) in which both players are allowed to use continuous as well as impulse controls. We are given $y_x(t):=y_x\bigl(t;\theta_1(.),u,\theta_2(.),v\bigr)$ the solution of (\ref{system}) which characterizes the state of the game at time $t$ with initial state $y_x(0^-)=x$ at $t=0^-$. The evolution of the system, described by the mapping $y_x:\mathbb{R}^+\rightarrow\mathbb{R}^n$, is provided by a deterministic model $\dot y_x(t)=b\bigl(y_x(t);\theta_1(t),\theta_2(t)\bigr)$, where $b$ is a function from $\mathbb{R}^n\times A\times B$ to $\mathbb{R}^n$ and $\theta_1(.)\in\Theta_1$ $\bigl(\text{resp.}\;\theta_2(.)\in\Theta_2\bigr)$ is the continuous control for $player-\xi$ (resp. $player-\eta$) defined in $\Theta_1$ (resp. $\Theta_2$) the space of measurable functions from $\mathbb{R}^+$ to $A$ (resp. $B$). The state undergoes impulses (jumps) $\xi_m$ and $\eta_k$, at certain impulse stopping times $\tau_m$ and $\rho_k$, respectively, that is:
\begin{equation*}
\left\{
\begin{aligned}
y_x(\tau_m^+)&=y_x(\tau_m^-)+g_\xi\bigl(y_x(\tau_m^-),\xi_m\bigr)\prod_{k\geq 0}\ind_{\{\tau_m\neq\rho_k\}},\;\tau_m\geq 0,\;\xi_m\neq 0;\\
y_x(\rho_k^+)&=y_x(\rho_k^-)+g_\eta\bigl(y_x(\rho_k^-),\eta_k\bigr),\;\rho_k\geq 0,\;\eta_k\neq 0,
\end{aligned}
\right.
\end{equation*}
where the impulse time sequences $\{\tau_m\}_{m\geq 0}$ and $\{\rho_k\}_{k\geq 0}$ are two non-decreasing sequences of $[0,\infty]$ such that $\tau_m,\rho_k\rightarrow+\infty$ when $m,k\rightarrow+\infty$, the impulse value sequences $\{\xi_m\}_{m\geq 0}$ and $\{\eta_k\}_{k\geq 0}$ are two sequences of elements of convex cones $U\subset\mathbb{R}^p$ and $V\subset\mathbb{R}^q$, respectively, and the form of impulses is of general term, that is it depends on the nonlinear function $g_\xi$ for $player-\xi$ and $g_\eta$ for $player-\eta$.
\par We call $\Theta_1\times\mathcal{U}$ and $\Theta_2\times\mathcal{V}$ the spaces of continuous and impulse controls $\bigl(\theta_1(.), u\bigr)$ and $\bigl(\theta_2(.), v\bigr)$ for $player-\xi$ and $player-\eta$, respectively. We denote $u:=(\tau_m,\xi_m)_{m\in\mathbb{N}}\in\mathcal{U}$ and $v:=(\rho_k,\eta_k)_{k\in\mathbb{N}}\in\mathcal{V}$ the two impulse controls such that $\xi_m\in U$ and $\eta_k\in V$. Thus $\bigl(\theta_1(.), u\bigr)\in\Theta_1\times\mathcal{U}$ and $\bigl(\theta_2(.), v\bigr)\in\Theta_2\times\mathcal{V}$ denote, respectively, the involved continuous and impulse controls for $player-\xi$ and $player-\eta$. For any initial state $x$ the controls $\bigl(\theta_1(.),u\bigr)$ and $\bigl(\theta_2(.),v\bigr)$ generate a trajectory $y_x(.)$ solution of the system (\ref{system}). Thus the state $y_x(.)$ is driven by two continuous and impulse controls, $\bigl(\theta_1(.),u:=(\tau_m,\xi_m)_{m\geq 0}\bigr)$ control of $player-\xi$ and $\bigl(\theta_2(.),v:=(\rho_k,\eta_k)_{k\geq 0}\bigr)$ control of $player-\eta$. The infinite product $\prod_{k\geq 0}\ind_{\{\tau_m\neq\rho_k\}}$ signifies that when the two players act together on the system at the same time, only the action of $player-\eta$ is tacking into account.
\par We are also given the following equation (\ref{payOff}), a gain (resp. cost) functional for $player-\xi$ (resp. $player-\eta$) as already mentioned in the introduction set:
\begin{equation}\label{payOff}\tag{$J$}
\begin{aligned}
J\bigl(x;\theta_1(.),u,\theta_2(.),v\bigr):=&\int_0^\infty f\bigl(y_x(t);\theta_1(t),\theta_2(t)\bigr)\exp(-\lambda t)dt\\
&-\sum_{m\geq 0}c\bigl(y_x(\tau_m^-),\xi_m\bigr)\exp(-\lambda\tau_m)\prod_{k\geq 0}\ind_{\{\tau_m\neq\rho_k\}}\\
&+\sum_{k\geq 0}\chi\bigl(y_x(\rho_k^-),\eta_k\bigr)\exp(-\lambda\rho_k),
\end{aligned}
\end{equation}
where $\bigl(\theta_1(.),u\bigr)\in\Theta_1\times\mathcal{U}$ and $\bigl(\theta_2(.),v\bigr)\in\Theta_2\times\mathcal{V}$ being the continuous and impulse controls. The functional $J$ represents the criterion which the $player-\xi$ wants to maximize and the $player-\eta$ wants to minimize. In the other words, $-J$ is the cost the $player-\eta$ has to pay, so the sum of the costs of the two players is null, which explains the name \textit{zero-sum}. We mention that $c$ and $\chi$ are the cost functions for $player-\xi$ and $player-\eta$, respectively, $f$ is the running gain and $\lambda>0$ the discount factor.
\par We assume that one player knows just the current and past choices of the control made by his opponent, thus we are given an information pattern for the two players prescribing that each of them choose his/her own control at each instant of time without knowing the future choices of the opponent. This is made rigorous by introducing, hereafter in the sense of Elliott-Kalton \cite{EK72,EK74}, the notion of \textit{non-anticipative strategy}.
\begin{definition}{(Non-Anticipative Strategy)}
A strategy for the $player-\xi$ is a map $\alpha:\Theta_2\times\mathcal{V}\rightarrow\Theta_1\times\mathcal{U}$; it is non-anticipative, if, for any $\theta_2^1(.),\theta_2^2(.)\in\Theta_2$, $v_1,v_2\in\mathcal{V}$ and $t>0$, $\theta_2^1(s)=\theta_2^2(s)$ and $v_1\equiv v_2$ on $[0,t]$ implies $\alpha\bigl(\theta_2^1(s),v_1\bigr)\equiv\alpha\bigl(\theta_2^2(s),v_2\bigr)$ for all $s\leq t$. We denote with $\mathcal{A}$ the set of all non-anticipative strategies $\alpha$ for $player-\xi$.\\
Similarly, the set of all non-anticipative strategies $\beta$ for $player-\eta$ is
\begin{equation*}
\begin{aligned}
\mathcal{B}:=\Bigl\{\beta:\Theta_1\times\mathcal{U}\rightarrow\Theta_2\times\mathcal{V}:\;&\theta_1^1(s)=\theta_1^2(s)\;\text{and}\;u_1\equiv u_2\;\text{on}\;[0,t]\;\text{for all}\;\theta_1^1(.),\theta_1^2(.)\in\Theta_1,\;u_1,u_2\in\mathcal{U},\\
&\;t>0\;\text{and}\;s\leq t;\text{implies}\;\beta\bigl(\theta_1^1(s),u_1\bigr)\equiv\beta\bigl(\theta_1^2(s),u_2\bigr)\Bigr\}.
\end{aligned}
\end{equation*}\qed
\end{definition}
Therefore we can define the lower and the upper value functions for the differential game control problem studied.
\begin{definition}
The lower value function of the game with the gain/cost functional $J:\mathbb{R}^n\times\Theta_1\times\mathcal{U}\times\Theta_2\times\mathcal{V}\rightarrow\mathbb{R}$ is
\begin{equation}\tag{$V^-$}\label{LowerValue}
V^-(x):=\inf_{\beta\in\mathcal{B}}\sup_{\bigl(\theta_1(.),u\bigr)\in\Theta_1\times\mathcal{U}}J\Bigl(x;\theta_1(.),u,\beta\bigl(\theta_1(.),u\bigr)\Bigr)
\end{equation}
and the upper value function is
\begin{equation}\tag{$V^+$}\label{UpperValue}
V^+(x):=\sup_{\alpha\in\mathcal{A}}\inf_{\bigl(\theta_2(.),v\bigr)\in\Theta_2\times\mathcal{V}}J\Bigl(x;\alpha\bigl(\theta_2(.),v\bigr),\theta_2(.),v\Bigr).
\end{equation}
\end{definition}
If $V^-(x)=V^+(x)$ we say that the game with initial point $x\in\mathbb{R}^n$ has a value, and we denote the value function of the game
\begin{equation}\tag{$V$}\label{value}
V(x):=V^-(x)=V^+(x).
\end{equation}\qed
\begin{remark}
Note that the inequality $V^-(x)\leq V^+(x)$ for all $x\in\mathbb{R}^n$, which would justify the terms \textit{lower} and \textit{upper}, is not obvious at first glance. Since in the definition of $V^-$ the inf is taken over non-anticipative strategies whereas in the definition of $V^+$ it is taken over controls, and similarly the sup is taken over different sets in the two definitions, then the inequality $V^+(x)\leq V^-(x)$ is false in general. We prove, in a rather indirect way by using the associated lower and upper HJBI QVIs (\ref{LowerQVI}) and (\ref{UpperQVI}), that the infinite horizon, zero-sum, deterministic differential game involving continuous and impulse controls studied in the present paper has a value. \qed
\end{remark}

\par Apart from the mathematical interest in its own right, the deterministic differential games enjoy a wide range of applications in various fields of engineering. We introduce, for the zero-sum games studied here, the following example:
\begin{example}{(Dynamic Portfolio Optimization)}\label{Example.1}
A typical example that provides an interesting framework of the theory of infinite horizon, zero-sum, deterministic games involving continuous and impulse controls, developed in the present paper, is the dynamic portfolio optimization problem described by the system (\ref{Spi}) below, where the market ($player-\xi$) is playing against the investor ($player-\eta$) and wishes to minimize his discounted terminal utility (i.e., maximize his cost defined by (\ref{Jpi}) below). We give the following illustrative dynamical system:
\begin{equation}\label{Spi}\tag{S$^\pi$}
\left\{
\begin{aligned}
y^\pi_x(0^-)&=x\in\mathbb{R}^n\;\text{(initial wealth)};\\
\dot y^\pi_x(t)&=b^\pi\bigl(y^\pi_x(t);\theta^\pi_1(t),\theta^\pi_2(t)\bigr),\;t\neq\tau_m,\;t\neq\rho_k;\\
y^\pi_x(\tau_m^+)&=y^\pi_x(\tau_m^-)+g_\xi^\pi\bigl(y^\pi_x(\tau_m^-),\xi_m\bigr)\prod_{k\geq 0}\ind_{\{\tau_m\neq\rho_k\}},\;\tau_m\geq 0,\;\xi_m\neq 0;\\
y^\pi_x(\rho_k^+)&=y^\pi_x(\rho_k^-)+g_\eta^\pi\bigl(y^\pi_x(\rho_k^-),\eta_k\bigr),\;\rho_k\geq 0,\;\eta_k\neq 0,
\end{aligned}
\right.
\end{equation}
where $x$ denotes the initial value of the investor's portfolio $\pi$, $y^\pi_x(t)$ represents the portfolio value (investor's wealth) at time $t$ controlled by:
\begin{enumerate}[i.]
\item Two continuous controls $\theta^\pi_1(.)$ and $\theta^\pi_2(.)$ which represent, respectively, the market's and the investor's instantaneous portfolio compositions. Thus, for example, $\theta^\pi_2(t)$ corresponds to the vector of number of units of stocks in the investor's portfolio at the instant $t$;
\item Two impulse controls $(\tau_m,\xi_m)_{m\in\mathbb{N}}$ and $(\rho_k,\eta_k)_{k\in\mathbb{N}}$ which describe, respectively, a new market's and investor's portfolio compositions at some jump instants. That is whenever the continuous controls $\theta_1^\pi(.)$ and $\theta_2^\pi(.)$ don't perform, the corresponding player uses a new optimal portfolio composition determined, respectively, at some impulse instants $\tau_m$ and $\rho_k$ with the impulse values $\xi_m$ and $\eta_k$, respectively.
\end{enumerate}
We denote $\bigl(\theta^\pi_1(.),u^\pi:=(\tau_m,\xi_m)_{m\in\mathbb{N}}\bigr)$ and $\bigl(\theta^\pi_2(.),v^\pi:=(\rho_k,\eta_k)_{k\in\mathbb{N}}\bigr)$ the continuous and impulse controls for the market and the investor, respectively, and we assume that the investor reacts immediately to the market whereas the market is not so quick in reacting to the investor's moves, i.e., the investor's action comes first whenever the impulse times for the two players coincide. Moreover, we assume that the investor does not consume wealth in the process of investing but is only interested to maximize his discounted terminal utility, that is, minimizing the following cost functional:
\begin{equation}\label{Jpi}\tag{$J^\pi$}
\begin{aligned}
J^\pi\bigl(x;\theta^\pi_1(.),u^\pi,\theta^\pi_2(.),v^\pi\bigr):=&\int_0^\infty f^\pi\bigl(y^\pi_x(t);\theta^\pi_1(t),\theta^\pi_2(t)\bigr)\exp(-\lambda t)dt\\
&-\sum_{m\geq 0}c^\pi\bigl(y^\pi_x(\tau_m^-),\xi_m\bigr)\exp(-\lambda\tau_m)\prod_{k\geq 0}\ind_{\{\tau_m\neq\rho_k\}}\\
&+\sum_{k\geq 0}\chi^\pi\bigl(y^\pi_x(\rho_k^-),\eta_k\bigr)\exp(-\lambda\rho_k),
\end{aligned}
\end{equation}
with the following components:
\begin{enumerate}[i.]
\item The running cost of integral type giving by the investor's stokes holding cost function $f^\pi$;
\item The maximizer's (market) $\bigl(\text{resp. minimizer's (investor)}\bigr)$ cost function $c^\pi$ $\bigl(\text{resp.}\;\chi^\pi\bigr)$ that corresponds to the cost of selling/buying stokes at impulse instants $\tau_m$ $(\text{resp.}\;\rho_k)$.
\end{enumerate}
The market moves according to the continuous control $\theta^\pi_1(.)$, creates jumps at impulse instants $\tau_m$ and tries to minimize the discounted terminal utility of the investor, that is, maximizing the cost functional (\ref{Jpi}), whereas the investor creates jumps at impulse instants $\rho_k$ and uses continuous control $\theta^\pi_2(.)$, obviously, tries to minimize the cost functional (\ref{Jpi}). We make the assumption that the flow of funds is between the investor and the market which makes our zero-sum game framework.\qed
\end{example}
Because of the advantage giving to the player using strategies, it is reasonable to believe that any more fair game has an outcome between $V^-(x)$ and $V^+(x)$. For this reason it is interesting to give the assumptions below ensuring the existence of a value, that is, the equality $V^-(x)=V^+(x)$ holds true for all $x\in\mathbb{R}^n$.

\subsection{Assumptions}\label{Sect.2.2}

In all the paper, we assume that $n$, $p$, $q$ and $l$ are some fixed positive integers, $k,m\in\mathbb{N}$, $T\in(0,+\infty]$, and we let the discount factor $\lambda$ be a fixed positive real. We denote, for $a,b\in\mathbb{R}$, $a\wedge b:=\min\{a,b\}$, by $|.|$ and $\|.\|$ the standard Euclidean vector norm in $\mathbb{R}$ and $\mathbb{R}^n$, respectively, and by $\|.\|_\infty$ the infinite norm in the space of bounded and continuous functions on $\mathbb{R}$ or $\mathbb{R}^n$.
\par We investigate the lower value (\ref{LowerValue}), the upper value (\ref{UpperValue}) of the differential game and the related HJBI QVIs (\ref{LowerQVI}) and (\ref{UpperQVI}) under the following less restrictive assumptions on the dynamics $b$, $g_\xi$, $g_\eta$, the running gain $f$ and the cost functions $c$ and $\chi$, where $\theta_1(.)\in\Theta_1$ and $\theta_2(.)\in\Theta_2$:
\begin{itemize}
\item[\textbf{[$\textbf{H}_b$]}] \textbf{(Dynamic)} We assume that the function $b$ is from $\mathbb{R}^n\times A\times B$, $\mathbb{R}^n-\text{valued}$, satisfies for some real constant $C_b>0$, all $x,y\in\mathbb{R}^n$ and all $t\geq 0$ the following:
$$\Bigl\|b\bigl(x;\theta_1(t),\theta_2(t)\bigr)-b\bigl(y;\theta_1(t),\theta_2(t)\bigr)\Bigr\|\leq C_b\|x-y\|,$$
and belongs to $C^0(\mathbb{R}^n\times A\times B)$, i.e., bounded and continuous in $\mathbb{R}^n\times A\times B$;
\item[\textbf{[$\textbf{H}_{g}$]}] \textbf{(Impulses Form)} We let the functions $g_\xi:(x,\xi)\in\mathbb{R}^n\times\mathbb{R}^p\rightarrow g_\xi(x,\xi)\in\mathbb{R}^n$ and $g_\eta:(x,\eta)\in\mathbb{R}^n\times\mathbb{R}^q\rightarrow g_\eta(x,\eta)\in\mathbb{R}^n$ be Lipschitz continuous with respect to $x$, uniformly in $\xi$ and $\eta$, respectively, with real constant $C_{g_\xi}>0$ and $C_{g_\eta}>0$, respectively;
\item[\textbf{[$\textbf{H}_f$]}] \textbf{(Running Gain)} We assume that the function $f$ is from $\mathbb{R}^n\times A\times B$, $\mathbb{R}-\text{valued}$, satisfies for some real constant $C_f>0$, all $x,y\in\mathbb{R}^n$ and all $t\geq 0$ the following:
$$\Bigl|f\bigl(x;\theta_1(t),\theta_2(t)\bigr)-f\bigl(y;\theta_1(t),\theta_2(t)\bigr)\Bigr|\leq C_f\|x-y\|,$$
and belongs to $C^0(\mathbb{R}^n\times A\times B)$;
\item[\textbf{[$\textbf{H}_{c,\chi}$]}] \textbf{(Impulses Cost)} The impulse cost functions $c:(x,\xi)\in\mathbb{R}^n\times U\rightarrow c(x,\xi)\in\mathbb{R}_{+}^{*}$ and $\chi:(x,\eta)\in\mathbb{R}^n\times V\rightarrow\chi(x,\eta)\in\mathbb{R}_{+}^{*}$ are from $\mathbb{R}^n$ and two convex cones $U\subset\mathbb{R}^p$ and $V\subset\mathbb{R}^q$, respectively, into $\mathbb{R}_{+}^{*}$, satisfy
\begin{equation}
\inf_{x\in\mathbb{R}^n,\;\xi\in U}c(x,\xi)>0\;\text{and}\;\inf_{x\in\mathbb{R}^n,\;\eta\in V}\chi(x,\eta)>0.
\end{equation}
and are Lipschitz continuous with respect to $x$, uniformly in $\xi$ and $\eta$, respectively, with real constants $C_{c}>0$ and $C_{\chi}>0$, respectively. Moreover, for all $x\in\mathbb{R}^n$, $\xi_1,\xi_2\in U$ and $\eta_1,\eta_2\in V$, we let the impulse costs satisfy
\begin{equation*}
\left\{
\begin{array}{ll}
c(x,\xi_1+\xi_2)&\leq c(x,\xi_1)+c(x,\xi_2);\\
\chi(x,\eta_1+\eta_2)&\leq\chi(x,\eta_1)+\chi(x,\eta_2),
\end{array}
\right.
\end{equation*}
that is multiple impulses occurring at the same time are sub-optimal.
\end{itemize}
\begin{remark}
It follows, regarding Assumption $\textbf{H}_{b}$, that there exists a unique global solution $y_{x}(.)$ to the above dynamical system (\ref{system}), while the Assumptions $\textbf{H}_{g}$, $\textbf{H}_{f}$ and $\textbf{H}_{c,\chi}$ provide the classical framework for the study, in the viscosity solutions framework, of the infinite horizon, zero-sum, deterministic continuous and impulse controls differential games considered. \qed
\end{remark}

\subsection{Classic Results}

\subsubsection{Dynamic Programming Principle}

Now we will prove the DPP property, initiated in the 1950s by Bellman \cite{B57}, in the following theorem, meaning that an optimal control viewed from today will remain optimal when viewed from tomorrow and stands for a most commonly used approach in solving optimal control problems:
\begin{theorem}{(Dynamic Programming Principle)}
Assume $\textbf{H}_b$, $\textbf{H}_g$, $\textbf{H}_f$ and $\textbf{H}_{c,\chi}$. Then for all $x\in\mathbb{R}^n$ and $T>0$, the lower value function (\ref{LowerValue}) and the upper value function (\ref{UpperValue}) satisfy, respectively,
\begin{equation}\label{equation_2}
\begin{aligned}
V^-(x)=\inf_{\beta\in\mathcal{B}}\sup_{\bigl(\theta_1(.),u\bigr)\in\Theta_1\times\mathcal{U}}\biggl\{&\int_0^T f\bigl(y_x(t);\theta_1(t),\theta_2(t)\bigr)\exp(-\lambda t)dt\\
&-\sum_{m\geq 0}c\bigl(y_x(\tau_m^-),\xi_m\bigr)\exp(-\lambda\tau_m)\ind_{\{\tau_m<T\}}\prod_{k\geq 0}\ind_{\{\tau_m\neq\rho_k\}}\\
&+\sum_{k\geq 0}\chi\bigl(y_x(\rho_k^-),\eta_k\bigr)\exp(-\lambda\rho_k)\ind_{\{\rho_k<T\}}\\
&+V^-\Bigl(y_x\bigl(T;\theta_1(.),u,\beta\bigl(\theta_1(.),u\bigr)\bigr)\Bigr)\exp(-\lambda T)\biggr\},
\end{aligned}
\end{equation}
where $\beta\bigl(\theta_1(.),u:=(\tau_m,\xi_m)_{m\in\mathbb{N}}\bigr):=\bigl(\theta_2(.),v:=(\rho_k,\eta_k)_{k\in\mathbb{N}}\bigr)$, and
\begin{equation*}
\begin{aligned}
V^+(x)=\sup_{\alpha\in\mathcal{A}}\inf_{\bigl(\theta_2(.),v\bigr)\in\Theta_2\times\mathcal{V}}\biggl\{&\int_0^T f\bigl(y_x(t);\theta_1(t),\theta_2(t)\bigr)\exp(-\lambda t)dt\\
&-\sum_{m\geq 0}c\bigl(y_x(\tau_m^-),\xi_m\bigr)\exp(-\lambda\tau_m)\ind_{\{\tau_m<T\}}\prod_{k\geq 0}\ind_{\{\tau_m\neq\rho_k\}}\\
&+\sum_{k\geq 0}\chi\bigl(y_x(\rho_k^-),\eta_k\bigr)\exp(-\lambda\rho_k)\ind_{\{\rho_k<T\}}\\
&+V^+\Bigl(y_x\bigl(T;\alpha\bigl(\theta_2(.),v\bigr),\theta_2(.),v\bigr)\Bigr)\exp(-\lambda T)\biggr\},
\end{aligned}
\end{equation*}
where $\alpha\bigl(\theta_2(.),v:=(\rho_k,\eta_k)_{k\in\mathbb{N}}\bigr):=\bigl(\theta_1(.),u:=(\tau_m,\xi_m)_{m\in\mathbb{N}}\bigr)$.
\end{theorem}
\begin{proof}
This proof is an adaptation of the results in chapter VIII of reference \cite{BC97}. We give only the proof for the lower value (\ref*{LowerValue}), similarly for the upper value (\ref*{UpperValue}). Let $T>0$, fix $\varepsilon>0$ and denote by $W_T(x)$ the right-hand side of (\ref{equation_2}). We first prove that $V^-(x)\leq W_T(x)$. For any $z\in\mathbb{R}^n$ we pick a non-anticipative strategy $\beta_z\in\mathcal{B}$ such that
\begin{equation}\label{equation_3}
V^-(z)\geq \sup_{\bigl(\theta_1(.),u\bigr)\in\Theta_1\times\mathcal{U}}J\Bigl(z;\theta_1(.),u,\beta_z\bigl(\theta_1(.),u\bigr)\Bigr)-\varepsilon,
\end{equation}
then we choose $\overline{\beta}\in\mathcal{B}$ a non-anticipative strategy for $player-\eta$ that satisfies, for $u:=(\tau_m,\xi_m)_{m\in\mathbb{N}}$, the following inequality:
\begin{equation}\label{equation_4}
\begin{aligned}
W_T(x)\geq& \sup_{\bigl(\theta_1(.),u\bigr)\in\Theta_1\times\mathcal{U}}\biggl\{\int_0^T f\bigl(y_x(t);\theta_1(t),\overline{\theta_2}(t)\bigr)\exp(-\lambda t)dt\\
&-\sum_{m\geq 0}c\bigl(y_x(\tau_m^-),\xi_m\bigr)\exp(-\lambda\tau_m)\ind_{\{\tau_m<T\}}\prod_{k\geq 0}\ind_{\{\tau_m\neq\overline{\rho_k}\}}\\
&+\sum_{k\geq 0}\chi\bigl(y_x(\overline{\rho_k}^-),\overline{\eta_k}\bigr)\exp(-\lambda\overline{\rho_k})\ind_{\{\rho_k<T\}}\\
&+V^-\Bigl(y_x\bigl(T;\theta_1(.),u,\overline{\beta}(\theta_1(.),u)\bigr)\Bigr)\exp(-\lambda T)\biggr\}-\varepsilon,
\end{aligned}
\end{equation}
where $$\overline{\beta}\bigl(\theta_1(.),u\bigr):=\bigl(\overline{\theta_2}(.),(\overline{\rho_k},\overline{\eta_k})_{k\in\mathbb{N}}\bigr).$$
Next, we define $\beta\in\mathcal{B}$, a non-anticipative strategy for $player-\eta$, as follows:
\begin{equation*}
\beta\bigl(\theta_1(t),u\bigr):=\left\{
\begin{aligned}
&\overline{\beta}\bigl(\theta_1(t),u\bigr),\;t\leq T;\\
&\beta_z\bigl(\theta_1^z(t-T),u^z\bigr),\;t>T,
\end{aligned}
\right.
\end{equation*}
where $z:=y_x\Bigl(T;\theta_1(.),u,\overline{\beta}\bigl(\theta_1(.),u\bigr)\Bigr)$, $\theta_1^z(.)\equiv\theta_1(.+T)$ and $u^z:=(\tau_m^z,\xi_m^z)_{m\in\mathbb{N}}$ with $\tau_m^z\geq T$. Since we have for all $t>0$, $$y_x\Bigl(t+T;\theta_1(.),u,\beta\bigl(\theta_1(.),u\bigr)\Bigr)=y_z\Bigl(t;\theta_1^z(.),u^z,\beta_z\bigl(\theta_1^z(.),u^z\bigr)\Bigr),$$
then by the change of variables $s=t+T$ we get
\begin{equation*}
\begin{aligned}
J\Bigl(z;\theta_1^z(.),u^z,\beta_z\bigl(\theta_1^z(.),u^z\bigr)\Bigr)=&\int_T^{+\infty} f\bigl(y_x(s);\theta_1(s),\theta_2(s)\bigr)\exp\bigl(-\lambda (s-T)\bigr)ds\\
&-\sum_{m\geq 0}c\bigl(y_x(\tau_m^-),\xi_m\bigr)\exp(-\lambda\tau_m)\ind_{\{\tau_m\geq T\}}\prod_{k\geq 0}\ind_{\{\tau_m\neq\rho_k\}}\\
&+\sum_{k\geq 0}\chi\bigl(y_x(\rho_k^-),\eta_k\bigr)\exp(-\lambda\rho_k)\ind_{\{\rho_k\geq T\}},
\end{aligned}
\end{equation*}
where $$\beta\bigl(\theta_1(.),u\bigr):=\bigl(\theta_2(.),(\rho_k,\eta_k)_{k\in\mathbb{N}}\bigr).$$
Then by (\ref{equation_3}) and (\ref{equation_4}) we deduce
\begin{equation*}
\begin{aligned}
W_T(x)\geq& \sup_{\bigl(\theta_1(.),u\bigr)\in\Theta_1\times\mathcal{U}}\biggl\{\int_0^{+\infty} f\bigl(y_x(t);\theta_1(t),\theta_2(t)\bigr)\exp(-\lambda t)dt\\
&-\sum_{m\geq 0}c\bigl(y_x(\tau_m^-),\xi_m\bigr)\exp(-\lambda\tau_m)\prod_{k\geq 0}\ind_{\{\tau_m\neq\rho_k\}}\\
&+\sum_{k\geq 0}\chi\bigl(y_x(\rho_k^-),\eta_k\bigr)\exp(-\lambda\rho_k)\biggr\}-2\varepsilon\\
&\geq V^-(x)-2\varepsilon,
\end{aligned}
\end{equation*}
thus, since $\varepsilon$ is arbitrary, we get the desired inequality.\\
We next prove that $W_T(x)\leq V^-(x)$. For any $z\in\mathbb{R}^n$ we pick the non-anticipative strategy $\beta_z\in\mathcal{B}$ for $player-\eta$ which satisfies the inequality (\ref{equation_3}). We pick $\bigl(\overline{\theta_1}(.),\overline{u}:=(\overline{\tau_m},\overline{\xi_m})_{m\in\mathbb{N}}\bigr)\in\Theta_1\times\mathcal{U}$, the continuous and impulse controls for $player-\xi$ that satisfy the following:
\begin{equation}\label{equation_5}
\begin{aligned}
W_T(x)\leq&\int_0^T f\bigl(y_x(t);\overline{\theta_1}(t),\theta_2^z(t)\bigr)\exp(-\lambda t)dt\\
&-\sum_{m\geq 0}c\bigl(y_x(\overline{\tau_m}^-),\overline{\xi_m}\bigr)\exp(-\lambda\overline{\tau_m})\ind_{\{\overline{\tau_m}<T\}}\prod_{k\geq 0}\ind_{\{\overline{\tau_m}\neq\rho_k^z\}}\\
&+\sum_{k\geq 0}\chi\bigl(y_x({\rho_k^z}^-),\eta_k^z\bigr)\exp(-\lambda\rho_k^z)\ind_{\{{\rho_k^z}<T\}}\\
&+V^-\Bigl(y_x\bigl(T;\overline{\theta_1}(.),\overline u,\beta_z\bigl(\overline{\theta_1}(.),\overline u\bigr)\bigr)\Bigr)\exp(-\lambda T)+\varepsilon,
\end{aligned}
\end{equation}
where $$\beta_z\bigl(\overline{\theta_1}(.),(\overline{\tau_m},\overline{\xi_m})_{m\in\mathbb{N}}\bigr):=\bigl(\theta_2^z(.),(\rho_k^z,\eta_k^z)_{k\in\mathbb{N}}\bigr).$$
For any $\theta_1(.)\in\Theta_1$ and $u\in\mathcal{U}$, we define the continuous control $\tilde{\theta_1}(.)\in\Theta_1$ for $player-\xi$ as follows:
\begin{equation}\label{equation_6}
\bigl(\tilde{\theta_1}(t),u\bigr):=\left\{
\begin{aligned}
&\bigl(\overline{\theta_1}(t),\overline u\bigr),\;t\leq T;\\
&\bigl(\theta_1(t-T),u\bigr),\;t>T,
\end{aligned}
\right.
\end{equation}
where $u:=(\tau_m,\xi_m)_{m\in\mathbb{N}}\in\mathcal{U}$ with $\tau_m\geq T$. Moreover, we define $\beta\in\mathcal{B}$ a non-anticipative strategy for $player-\eta$ as follows:
\begin{equation}\label{equation_7}
\beta\bigl(\theta_1(t),u\bigr):=\beta_z\bigl(\tilde{\theta_1}(t+T),u\bigr).
\end{equation}
Next, set
\begin{equation}\label{equation_8}
z_1:=y_x\Bigl(T;\overline{\theta_1}(.),\overline u,\beta_{z}\bigl(\overline{\theta_1}(.),\overline u\bigr)\Bigr),
\end{equation}
and choose $\theta_1(.)\in\Theta_1$ and $u\in\mathcal{U}$ such that
\begin{equation}\label{equation_9}
V^-(z_1)\leq J\Bigl(z_1;\theta_1(.),u,\beta\bigl(\theta_1(.),u\bigr)\Bigr)+\varepsilon.
\end{equation}
Observe that, by (\ref{equation_6}) and (\ref{equation_7}), we have
\begin{equation*}
y_x\Bigl(s;\tilde{\theta_1}(.),u,\beta_z\bigl(\tilde{\theta_1}(.),u\bigr)\Bigr)=\left\{
\begin{aligned}
&y_x\Bigl(s;\overline{\theta_1}(.),\overline u,\beta_z\bigl(\overline{\theta_1}(.),\overline u\bigr)\Bigr),\;s\leq T;\\
&y_{z_1}\Bigl(s-T;\theta_1(.),u,\beta\bigl(\theta_1(.),u\bigr)\Bigr),\;s>T,
\end{aligned}
\right.
\end{equation*}
so by the change of variable $s=t+T$ we deduce for $u:=(\tau_m,\xi_m)_{m\in\mathbb{N}}$ that
\begin{equation}\label{equation_10}
\begin{aligned}
J\Bigl(z_1;\theta_1(.),u,\beta\bigl(\theta_1(.),u\bigr)\Bigr)=&\int_T^{+\infty} f\bigl(y_x(s);\tilde{\theta_1}(s),\tilde{\theta_2}(s)\bigr)\exp\bigl(-\lambda (s-T)\bigr)ds\\
&-\sum_{m\geq 0}c\bigl(y_x(\tau_m^-),\xi_m\bigr)\exp(-\lambda\tau_m)\ind_{\{\tau_m\geq T\}}\prod_{k\geq 0}\ind_{\{\tau_m\neq\tilde{\rho}_k\}}\\
&+\sum_{k\geq 0}\chi\bigl(y_x(\tilde{\rho}_k^-),\tilde{\eta}_k\bigr)\exp(-\lambda\tilde{\rho}_k)\ind_{\{\tilde{\rho}_k\geq T\}},
\end{aligned}
\end{equation}
where $$\beta_z\bigl(\tilde{\theta_1}(.),u\bigr):=\bigl(\tilde{\theta_2}(.),(\tilde{\rho}_k,\tilde{\eta}_k)_{k\in\mathbb{N}}\bigr),$$
Now we use (\ref{equation_5}), (\ref{equation_6}), (\ref{equation_8}), (\ref{equation_9}) and (\ref{equation_10}) to get
\begin{equation*}
W_T(x)\leq J\Bigl(x;\tilde{\theta_1}(.),u,\beta_z\bigl(\tilde{\theta_1}(.),u\bigr)\Bigr)+2\varepsilon,
\end{equation*}
thus, from inequality (\ref{equation_3}), we deduce that $W_T(x)\leq V^-(x)+3\varepsilon$. Then, since $\varepsilon$ is arbitrary, we obtain the desired inequality.
\end{proof}

\subsubsection{Regularity of the Value Functions}
We prove hereafter some results concerning the boundedness and the regularity of the lower value (\ref{LowerValue}) and the upper value (\ref{UpperValue}). We start by an estimate on the trajectories. Let $x,z\in\mathbb{R}^n$ and denote $y_x\bigl(.;\theta_1(.),u,\theta_2(.),v\bigr)$ and $y_z\bigl(.;\theta_1(.),u,\theta_2(.),v\bigr)$ the two trajectories generated, respectively, from $x$ and $z$ by the continuous and impulse controls $\bigl(\theta_1(.),u:=(\tau_m,\xi_m)_{m\in\mathbb{N}}\bigr)\in\Theta_1\times\mathcal{U}$ and $\bigl(\theta_2(.),v:=(\rho_k,\eta_k)_{k\in\mathbb{N}}\bigr)\in\Theta_2\times\mathcal{V}$, we then have the following estimate:
\begin{proposition}\label{Proposition1}
Assume $\textbf{H}_b$ and $\textbf{H}_g$. We have, for all $x,z\in\mathbb{R}^n$ and $t\geq 0$, the usual estimate on the trajectories:
$$\Bigl\|y_x\bigl(t;\theta_1(.),u,\theta_2(.),v\bigr)-y_z\bigl(t;\theta_1(.),u,\theta_2(.),v\bigr)\Bigr\|\leq\exp(Ct)\|x-z\|,$$
where $C$ is a positive real constant depending on the number of impulses.
\end{proposition}
\begin{proof}
Let $\bigl(\theta_1(.),u\bigr)\in\Theta_1\times\mathcal{U}$ and $\bigl(\theta_2(.),v\bigr)\in\Theta_2\times\mathcal{V}$. By Gronwall's Lemma, using Assumption $\textbf{H}_b$ we have for all $t\in[0,\tau_0\wedge\rho_0]$,
\begin{equation*}
\Bigl\|y_x\bigl(t;\theta_1(.),u,\theta_2(.),v\bigr)-y_z\bigl(t;\theta_1(.),u,\theta_2(.),v\bigr)\Bigr\|\leq\exp(C_bt)\|x-z\|,
\end{equation*}
moreover, from Assumption $\textbf{H}_g$, we have
\begin{equation*}
\begin{aligned}
\Bigl\|y_x\bigl(\tau_0^+\wedge\rho_0^+;\theta_1(.),u,\theta_2(.),v\bigr)-&y_z\bigl(\tau_0^+\wedge\rho_0^+;\theta_1(.),u,\theta_2(.),v\bigr)\Bigr\|\\
&\leq\;(1+C_g)\exp\bigl(C_b(\tau_0\wedge\rho_0)\bigr)\|x-z\|,
\end{aligned}
\end{equation*}
where $$C_g=C_{g_\xi}\ind_{\{\tau_0<\rho_0\}}+C_{g_\eta}\ind_{\{\tau_0\geq\rho_0\}}.$$\\
Repeating inductively the same argument to get, for the impulse time $\tau_m\wedge\rho_m$, that
\begin{equation*}
\begin{aligned}
\Bigl\|y_x\bigl(\tau_m^+\wedge\rho_m^+;\theta_1(.),u,\theta_2(.),v\bigr)-&y_z\bigl(\tau_m^+\wedge\rho_m^+;\theta_1(.),u,\theta_2(.),v\bigr)\Bigr\|\\
&\leq\;(1+2\tilde{C_g})^{N}\exp\bigl(C_b(\tau_m\wedge\rho_m)\bigr)\|x-z\|,
\end{aligned}
\end{equation*}
where the constant $N$ depends on the number of impulses $m$ and given by
\begin{equation*}
N=m-\sum_{i=0}^{m}\ind_{\{\tau_i=\rho_i\}},\;\text{and}\;\tilde{C_g}=\max\bigl\{C_{g_\xi},C_{g_\eta}\bigr\}.
\end{equation*}
Thus, for all $t\in[\tau_m\wedge\rho_m,\tau_{m+1}\wedge\rho_{m+1}]$, we get the existence of a constant $C>0$, depending on the number of impulses, such that
$$\Bigl\|y_x\bigl(t;\theta_1(.),u,\theta_2(.),v\bigr)-y_z\bigl(t;\theta_1(.),u,\theta_2(.),v\bigr)\Bigr\|\leq\exp(Ct)\|x-z\|,$$
this inequality remains true even when $t$ is greater than the last impulse time. The proof is then complete.
\end{proof}
We are going now to prove the following theorem:
\begin{theorem}
Assume $\textbf{H}_b$, $\textbf{H}_g$, $\textbf{H}_f$ and $\textbf{H}_{c,\chi}$. Then the lower value function (\ref{LowerValue}) and the upper value function (\ref{UpperValue}) are in the space of bounded and uniformly continuous functions in $\mathbb{R}^n$.
\end{theorem}
\begin{proof}
We give only the proof for the lower value (\ref*{LowerValue}), similarly for the upper value (\ref*{UpperValue}). We proceed for the proof in two steps:\\
\par \textbf{Step 1: Boundedness.} Let $x\in\mathbb{R}^n$ and $\beta\in\mathcal B$ be any non-anticipative strategy for $player-\eta$, we have
$$V^-(x)\leq \sup_{\bigl(\theta_1(.),u\bigr)\in\Theta_1\times\mathcal{U}}J\Bigl(x;\theta_1(.),u,\beta\bigl(\theta_1(.),u\bigr)\Bigr),$$
considering the set of non-anticipative strategies $\beta\bigl(\theta_1(.),u\bigr):=\bigl(\theta_2(.),(\rho_k,\eta_k)_{k\in\mathbb{N}}\bigr)$ where there is no impulse time, i.e., $\rho_0=+\infty$, for $u:=(\tau_m,\xi_m)_{m\in\mathbb{N}}$, we get
\begin{equation*}
V^-(x)\leq \sup_{\bigl(\theta_1(.),u\bigr)\in\Theta_1\times\mathcal{U}}\biggl\{\int_0^{+\infty} f\bigl(y_x(t);\theta_1(t),\theta_2(t)\bigr)\exp(-\lambda t)dt-\sum_{m\geq 0}c\bigl(y_x(\tau_m^-),\xi_m\bigr)\exp(-\lambda\tau_m)\biggr\}.
\end{equation*}
Fix now $\varepsilon>0$, then there exists a continuous and impulse control $\bigl(\overline{\theta_1}(.),\overline{u}:=(\overline{\tau_m},\overline{\xi_m})\bigr)\in\Theta_1\times\mathcal{U}$ for $player-\xi$ such that
\begin{equation*}
V^-(x)\leq \int_0^{+\infty} f\bigl(y_x(t);\overline{\theta_1}(t),\theta_2(t)\bigr)\exp(-\lambda t)dt-\sum_{m\geq 0}c\bigl(y_x(\overline{\tau_m}^-),\overline{\xi_m}\bigr)\exp(-\lambda\overline{\tau_m})+\varepsilon.
\end{equation*}
Since, from Assumptions $\textbf{H}_{c,\chi}$ and $\textbf{H}_f$, $c$ is a non negative function and $f$ is bounded and since the term $\sum c\exp$ is bounded for infinite number of impulses, then there exists $C>0$ such that
$$V^-(x)\leq C.$$
Similarly, wet let $\bigl(\theta_1(.),u:=(\tau_m,\xi_m)_{m\in\mathbb{N}}\bigr)\in\Theta_1\times\mathcal{U}$ be the set of continuous and impulse controls for $player-\xi$ for which there is no impulse time, i.e., $\tau_0=+\infty$ and which satisfies
\begin{equation*}
V^-(x)\geq \inf_{\beta\in\mathcal{B}}\biggl\{\int_0^{+\infty} f\bigl(y_x(t);\theta_1(t),\theta_2(t)\bigr)\exp(-\lambda t)dt+\sum_{k\geq 0}\chi\bigl(y_x(\rho_k^-),\eta_k\bigr)\exp(-\lambda\rho_k)\biggr\},
\end{equation*}
where $$\beta\bigl(\theta_1(.),u\bigr):=\bigl(\theta_2(.),v:=(\rho_k,\eta_k)_{k\in\mathbb{N}}\bigr).$$
There exists a non-anticipative strategy $\overline{\beta}\in\mathcal{B}$ for $player-\eta$ which provides a continuous and impulse control $\bigl(\overline{\theta_2}(.),\overline{v}:=(\overline{\rho_k},\overline{\eta_k})\bigr)\in\Theta_2\times\mathcal{V}$ that satisfies, for a fixed $\varepsilon>0$,
\begin{equation*}
V^-(x)\geq \int_0^{+\infty} f\bigl(y_x(t);\theta_1(t),\overline{\theta_2}(t)\bigr)\exp(-\lambda t)dt+\sum_{k\geq 0}\chi\bigl(y_x(\overline{\rho_k}^-),\overline{\eta_k}\bigr)\exp(-\lambda\overline{\rho_k})-\varepsilon.
\end{equation*}
Recall Assumptions $\textbf{H}_{c,\chi}$ and $\textbf{H}_f$, $\chi$ is a non negative function and $f$ is bounded and since the term $\sum\chi\exp$ is bounded for infinite number of impulses, then there exists $C^{'}>0$ such that
$$V^-(x)\geq -C^{'}.$$
Which finishes the proof.\\
\par \textbf{Step 2: Uniform Continuity.} We fix $x,z\in\mathbb{R}^n$, $\varepsilon>0$ and first pick $\overline{\beta}\in\mathcal{B}$ a non-anticipative strategy for $player-\eta$ that satisfies
$$V^-(z)\geq \sup_{\bigl(\theta_1(.),u\bigr)\in\Theta_1\times\mathcal{U}}J\Bigl(z;\theta_1(.),u,\overline{\beta}\bigl(\theta_1(.),u\bigr)\Bigr)-\frac{\varepsilon}{2},$$
then we pick $\bigl(\overline{\theta_1}(.),\overline{u}:=(\overline{\tau_m},\overline{\xi_m})_{m\in\mathbb{N}}\bigr)\in\Theta_1\times\mathcal{U}$, the continuous and impulse controls for $player-\xi$ that satisfies
\begin{equation*}
\begin{aligned}
V^-(x)&\leq\sup_{\bigl(\theta_1(.),u\bigr)\in\Theta_1\times\mathcal{U}}J\Bigl(x;\theta_1(.),u,\overline{\beta}\bigl(\theta_1(.),u\bigr)\Bigr)\\
&\leq J\Bigl(x;\overline{\theta_1}(.),\overline{u},\overline{\beta}\bigl(\overline{\theta_1}(.),\overline{u}\bigr)\Bigr)+\frac{\varepsilon}{2}.
\end{aligned}
\end{equation*}
Thus
$$V^-(x)-V^-(z)\leq J\Bigl(x;\overline{\theta_1}(.),\overline{u},\overline{\beta}\bigl(\overline{\theta_1}(.),\overline{u}\bigr)\Bigr)-J\Bigl(z;\overline{\theta_1}(.),\overline{u},\overline{\beta}\bigl(\overline{\theta_1}(.),\overline{u}\bigr)\Bigr)+\varepsilon.$$
It follows
\begin{equation*}
\begin{aligned}
V^-(x)-V^-(z)\leq&\int_0^{+\infty} \Bigl[f\bigl(y_x(t);\overline{\theta_1}(t),\overline{\theta_2}(t)\bigr)-f\bigl(y_z(t);\overline{\theta_1}(t),\overline{\theta_2}(t)\bigr)\Bigr]\exp(-\lambda t)dt\\
&-\sum_{m\geq 0}c\bigl(y_x(\overline{\tau_m}^-),\overline{\xi_m}\bigr)\exp(-\lambda\overline{\tau_m})\prod_{k\geq 0}\ind_{\{\overline{\tau_m}\neq\overline{\rho_k}\}}\\
&+\sum_{k\geq 0}\chi\bigl(y_x(\overline{\rho_k}^-),\overline{\eta_k}\bigr)\exp(-\lambda\overline{\rho_k})\\
&+\sum_{m\geq 0}c\bigl(y_z(\overline{\tau_m}^-),\overline{\xi_m}\bigr)\exp(-\lambda\overline{\tau_m})\prod_{k\geq 0}\ind_{\{\overline{\tau_m}\neq\overline{\rho_k}\}}\\
&-\sum_{k\geq 0}\chi\bigl(y_z(\overline{\rho_k}^-),\overline{\eta_k}\bigr)\exp(-\lambda\overline{\rho_k})+\varepsilon,
\end{aligned}
\end{equation*}
where $\overline{\beta}\bigl(\overline{\theta_1}(.),\overline{u}\bigr):=\bigl(\overline{\theta_2}(.),\overline{v}:=(\overline{\rho}_k,\overline{\eta}_k)_{k\in\mathbb{N}}\bigr)\in\Theta_2\times\mathcal{V}$. Then, using the definition of the lower value (\ref*{LowerValue}) given by the DPP property (\ref{equation_2}) for $T>0$, we obtain
\begin{equation*}
\begin{aligned}
V^-(x)-V^-(z)\leq&\int_0^T\Bigl[f\bigl(y_x(t);\overline{\theta_1}(t),\overline{\theta_2}(t)\bigr)-f\bigl(y_{z}(t);\overline{\theta_1}(t),\overline{\theta_2}(t)\bigr)\Bigr]\exp(-\lambda t)dt\\
&-\sum_{m\geq 0}c\bigl(y_x(\overline{\tau_m}^-),\overline{\xi_m}\bigr)\exp(-\lambda\overline{\tau_m})\ind_{\{\overline{\tau_m}<T\}}\prod_{k\geq 0}\ind_{\{\overline{\tau_m}\neq\overline{\rho_k}\}}\\
&+\sum_{k\geq 0}\chi\bigl(y_x(\overline{\rho_k}^-),\overline{\eta_k}\bigr)\exp(-\lambda\overline{\rho_k})\ind_{\{\overline{\rho_k}<T\}}\\
&+\sum_{m\geq 0}c\bigl(y_z(\overline{\tau_m}^-),\overline{\xi_m}\bigr)\exp(-\lambda\overline{\tau_m})\ind_{\{\overline{\tau_m}<T\}}\prod_{k\geq 0}\ind_{\{\overline{\tau_m}\neq\overline{\rho_k}\}}\\
&-\sum_{k\geq 0}\chi\bigl(y_z(\overline{\rho_k}^-),\overline{\eta_k}\bigr)\exp(-\lambda\overline{\rho_k})\ind_{\{\overline{\rho_k}<T\}}\\
&+V^-\Bigl(y_x\bigl(T;\overline{\theta_1}(.),\overline{u},\overline{\beta}\bigl(\overline{\theta_1}(.),\overline{u}\bigr)\bigr)\Bigr)\exp(-\lambda T)\\
&-V^-\Bigl(y_z\bigl(T;\overline{\theta_1}(.),\overline{u},\overline{\beta}\bigl(\overline{\theta_1}(.),\overline{u}\bigr)\bigr)\Bigr)\exp(-\lambda T)+\varepsilon.
\end{aligned}
\end{equation*}
Thus, by Assumption $\textbf{H}_f$ and Assumption $\textbf{H}_{c,\chi}$, we get
\begin{equation*}
\begin{aligned}
V^-(x)-V^-(z)\leq&\int_0^T C_f\bigl\|y_x(t)-y_{z}(t)\bigr\|\exp(-\lambda t)dt\\
&-\sum_{m\geq 0}C_c\bigl\|y_x(\overline{\tau_m}^-)-y_z(\overline{\tau_m}^-)\bigr\|\exp(-\lambda\overline{\tau_m})\ind_{\{\overline{\tau_m}<T\}}\prod_{k\geq 0}\ind_{\{\overline{\tau_m}\neq\overline{\rho_k}\}}\\
&+\sum_{k\geq 0}C_{\chi}\bigl\|y_x(\overline{\rho_k}^-)-y_z(\overline{\rho_k}^-)\bigr\|\exp(-\lambda\overline{\rho_k})\ind_{\{\overline{\rho_k}<T\}}\\
&+\Bigl|V^-\Bigl(y_x\bigl(T;\overline{\theta_1}(.),\overline{u},\overline{\beta}\bigl(\overline{\theta_1}(.),\overline{u}\bigr)\bigr)\Bigr)-V^-\Bigl(y_z\bigl(T;\overline{\theta_1}(.),\overline{u},\overline{\beta}\bigl(\overline{\theta_1}(.),\overline{u}\bigr)\bigr)\Bigr)\Bigr|\\
&\times\exp(-\lambda T)+\varepsilon.
\end{aligned}
\end{equation*}
By Proposition \ref{Proposition1} and the boundedness of $V^-$, we deduce that there exist some constants $C,C_v>0$ such that
\begin{equation}\label{equation_11}
\begin{aligned}
V^-(x)-V^-(z)\leq&\;C_f\|x-z\|\int_0^T \exp\bigl((C-\lambda)t\bigr)dt\\
&-C_c\|x-z\|\sum_{m\geq 0,\;\overline{\tau_m}<T}\exp\bigl((C-\lambda)\overline{\tau_m}\bigr)\prod_{k\geq 0}\ind_{\{\overline{\tau_m}\neq\overline{\rho_k}\}}\\
&+C_{\chi}\|x-z\|\sum_{k\geq 0,\;\overline{\rho_k}<T}\exp\bigl((C-\lambda)\overline{\rho_k}\bigr)\\
&+2C_v\exp(-\lambda T)+\varepsilon.
\end{aligned}
\end{equation}
Therefore, if $\lambda>C$, since in the right-hand side of (\ref{equation_11}) the sums are finite, there exists $0<K<+\infty$ such that
\begin{equation}\label{equation_12}
\begin{aligned}
V^-(x)-V^-(z)\leq&\;\frac{C_f}{C-\lambda}\| x-z\|\Bigl[\exp\bigl((C-\lambda)T\bigr)-1\Bigr]\\
&+K\|x-z\|+2C_v\exp(-\lambda T)+\varepsilon.
\end{aligned}
\end{equation}
Hence, since in the right-hand side of (\ref{equation_12}) the terms $\exp\bigl((C-\lambda)T\bigr)$ and $2C_v\exp(-\lambda T)$ are bounded, we deduce from the arbitrariness of $\varepsilon$ that there exists a constant $M>0$ such that
$$V^-(x)-V^-(z)\leq M\|x-z\|.$$
In the case where $\lambda<C$, we choose $T$ such that $\exp(-CT)=\|x-z\|^{1/2}$ with $\|x-z\|<1$. Hence, in the right-hand side of (\ref{equation_12}), the first term becomes
$$\frac{C_f}{C-\lambda}\|x-z\|^{1/2}\bigl(\exp(-\lambda T)-\|x-z\|^{1/2}\bigr),$$
and the term $2C_v\exp(-\lambda T)$ is bounded. We then deduce, from the arbitrariness of $\varepsilon$, the same last inequality. Finally, when $\lambda=C$, it suffice to let some $\hat{\lambda}<\lambda=C$, so we go back to (\ref{equation_11}) and we proceed, since $\exp\bigl((C-\lambda) T\bigr)<\exp\bigl((C-\hat{\lambda})T\bigr)$ and $\exp(-\lambda T)<\exp(-\hat{\lambda}T)$, as above with the case $\hat{\lambda}<C$.
Hence, from the fact that the role of $x$ and $z$ being symmetrical, one might deduce the uniformly continuity of the lower value ($V^-$) in $\mathbb{R}^n$.
\end{proof}
Next, we prove the following useful proposition:
\begin{proposition}\label{proposition2}
If the function $v$ is uniformly continuous in $\mathbb{R}^n$, so is the two functions $\mathcal{H}_{inf}^\chi v$ and $\mathcal{H}_{sup}^c v$.
\end{proposition}
\begin{proof}
We give the proof for $\mathcal{H}_{inf}^\chi v$, similarly for $\mathcal{H}_{sup}^c v$. Let $v$ be a uniformly continuous function, $x,z\in\mathbb{R}^n$ and choose $\varepsilon>0$ and $\eta_\varepsilon\in V$ such that we have
$$\mathcal{H}_{inf}^\chi v(z)+\varepsilon\geq v\bigl(z+g_\eta(z,\eta_\varepsilon)\bigr)+\chi(z,\eta_\varepsilon),$$
thus
$$\mathcal{H}_{inf}^\chi v(x)-\mathcal{H}_{inf}^\chi v(z)\leq v\bigl(x+g_\eta(x,\eta_\varepsilon)\bigr)+\chi(x,\eta_\varepsilon)-v\bigl(z+g_\eta(z,\eta_\varepsilon)\bigr)-\chi(z,\eta_\varepsilon)+\varepsilon.$$
It follows, from Assumption $\textbf{H}_{c,\chi}$, the existence of a constant $C>0$ such that
$$\mathcal{H}_{inf}^\chi v(x)-\mathcal{H}_{inf}^\chi v(z)\leq C\|x-z\|+\varepsilon,$$
since $x$ and $z$ play symmetrical roles, we deduce from the arbitrariness of $\varepsilon$ that $\mathcal{H}_{inf}^\chi v$ is uniformly continuous.
\end{proof}

\section{Viscosity Characterization of the Game}\label{Sect.3}

For the two-player, zero-sum, deterministic, continuous and impulse controls game studied the associated lower and upper HJBI QVIs are derived from the DPP and given, respectively, by the aforementioned equations (\ref{LowerQVI}) and (\ref{UpperQVI}).

\begin{remark}
In this paper, our aim is to show that the differential game considered has a value, and we proceed in two steps:
\begin{enumerate}[i.]
\item First, we study the existence of the solution in viscosity sense for both HJBI QVIs (\ref{LowerQVI}) and (\ref{UpperQVI}), i.e., we prove the fact that the lower value (\ref{LowerValue}) (resp. upper value (\ref{UpperValue})) is a viscosity solution of the lower HJBI QVI (\ref{LowerQVI}) (resp. upper HJBI QVI (\ref{UpperQVI}));
\item Next, we show that both lower HJBI QVI (\ref{LowerQVI}) and upper HJBI QVI (\ref{UpperQVI}) admit, under Isaacs' condition, the lower value (\ref{LowerValue}) and the upper value (\ref{UpperValue}) as unique solutions in the viscosity sense.
\end{enumerate}
Thus, the game admits the value function (\ref{value}), i.e., $V^-(x)=V^+(x)$ for all $x\in\mathbb{R}^n$. \qed
\end{remark}

\par In the rest of this section and in Sect. \ref{Sect.4}, we adopt the following definition of the viscosity solution:
\begin{definition}{(Viscosity Solution)}
Let $V:\mathbb{R}^n\rightarrow\mathbb{R}$ be a continuous function. $V$ is called:
\begin{enumerate}[i.]
\item A viscosity sub-solution of the lower HJBI QVI (\ref{LowerQVI}) (resp. upper HJBI QVI (\ref{UpperQVI})) if for any $\overline x\in\mathbb{R}^n$ and any function $\phi\in C^{1}(\mathbb{R}^n)$ such that $\overline x$ is a local maximum point of $V-\phi$, we have:
$$\min\biggl\{\max\Bigl[\lambda V(\overline x)+H^-\bigl(\overline x,D\phi(\overline x)\bigr),V(\overline x)-\mathcal{H}_{inf}^\chi V(\overline x)\Bigr],V(\overline x)-\mathcal{H}_{sup}^c V(\overline x)\biggr\}\leq 0$$
$$\biggl(\text{resp}.\;\max\biggl\{\min\Bigl[\lambda V(\overline x)+H^+\bigl(\overline x,D\phi(\overline x)\bigr),V(\overline x)-\mathcal{H}_{sup}^c V(\overline x)\Bigr],V(\overline x)-\mathcal{H}_{inf}^\chi V(\overline x)\biggr\}\leq 0\biggr);$$
\item A viscosity super-solution of the lower HJBI QVI (\ref{LowerQVI}) (resp. upper HJBI QVI (\ref{UpperQVI})) if for any $\underline x\in\mathbb{R}^n$ and any function $\phi\in C^{1}(\mathbb{R}^n)$ such that $\underline x$ is a local minimum point of $V-\phi$, we have:
$$\min\biggl\{\max\Bigl[\lambda V(\underline x)+H^-\bigl(\underline x,D\phi(\underline x)\bigr),V(\underline x)-\mathcal{H}_{inf}^\chi V(\underline x)\Bigr],V(\underline x)-\mathcal{H}_{sup}^c V(\underline x)\biggr\}\geq 0$$
$$\biggl(\text{resp}.\;\max\biggl\{\min\Bigl[\lambda V(\underline x)+H^+\bigl(\underline x,D\phi(\underline x)\bigr),V(\underline x)-\mathcal{H}_{sup}^c V(\underline x)\Bigr],V(\underline x)-\mathcal{H}_{inf}^\chi V(\underline x)\biggr\}\geq 0\biggr);$$
\item A viscosity solution of the lower HJBI QVI (\ref{LowerQVI}) (resp. upper HJBI QVI (\ref{UpperQVI})) if it is both a viscosity sub-solution and super-solution of lower HJBI QVI (\ref{LowerQVI}) (resp. upper HJBI QVI (\ref{UpperQVI})).\qed
\end{enumerate}
\end{definition}

\subsection{Useful Lemmas}

Next, we give the proof for the following lemmas which will be useful to prove the existence of the viscosity solution for the HJBI QVIs (\ref{LowerQVI}) and (\ref{UpperQVI}):
\begin{lemma}\label{lemma_1}
Assume $\textbf{H}_b$, $\textbf{H}_g$, $\textbf{H}_f$ and $\textbf{H}_{c,\chi}$. The lower value function (\ref{LowerValue}) satisfies for all $x\in\mathbb{R}^n$ the following properties:
\begin{enumerate}[i.]
\item $V^-(x)\leq\mathcal{H}_{inf}^\chi V^-(x)$;
\item If $V^-(x)<\mathcal{H}_{inf}^\chi V^-(x)$ then $V^-(x)\geq\mathcal{H}_{sup}^c V^-(x)$.
\end{enumerate}
The same properties hold true for the upper value function (\ref{UpperValue}).
\end{lemma}
\begin{proof}
We give only the proof for $V^-$, similarly for $V^+$. First, we let $x\in\mathbb{R}^n$, $\theta_1(.)\in\Theta_1$ and $u:=(\tau_m,\xi_m)_{m\in\mathbb{N}}\in\mathcal{U}$, then we consider, for $player-\eta$, the non-anticipative strategy $\beta\bigl(\theta_1(.),u\bigr):=\bigl(\theta_2(.),v:=(\rho_k,\eta_k)_{k\in\mathbb{N}}\bigr)\in\Theta_2\times\mathcal{V}$ where $\beta\in\mathcal{B}$. Next, choose $\beta^{'}\in\mathcal{B}$ such that $\beta^{'}(.,u):=\bigl(.,(0,\eta;\rho_1,\eta_1;\rho_2,\eta_2;...)\bigr)$, we then obtain $$V^-(x)\leq\sup_{\bigl(\theta_1(.),u\bigr)\in\Theta_1\times\mathcal{U}}J\Bigl(x;\theta_1(.),u,\beta^{'}\bigl(\theta_1(.),u\bigr)\Bigr),$$
thus
$$V^-(x)\leq\sup_{\bigl(\theta_1(.),u\bigr)\in\Theta_1\times\mathcal{U}}J\Bigl(x+g_\eta(x,\eta);\theta_1(.),u,\beta\bigl(\theta_1(.),u\bigr)\Bigr)+\chi(x,\eta),$$
from which we get
$$V^-(x)\leq V^-\bigl(x+g_\eta(x,\eta)\bigr)+\chi(x,\eta).$$
then the inequality (\textbf{$\textit{i}$}) follows from the arbitrariness of $\eta$.\\
Now let us assume that $V^-(x)<\mathcal{H}_{inf}^\chi V^-(x)$ for some $x\in\mathbb{R}^n$. From the DPP for the lower value $V^-$ (\ref{equation_2}), by taking $T=0$ with possible impulses at $0$, the sub-optimality of multiple impulses at the same time, and Assumption $\textbf{H}_{c,\chi}$ we get
\begin{equation*}
\begin{aligned}
V^{-}(x)=&\inf_{\underset{\rho_0\in\{0,+\infty\},\;\eta\in V}{\theta_2(.)\in\Theta_2}}\sup_{\underset{\tau_0\in\{0,+\infty\},\;\xi\in U}{\theta_1(.)\in\Theta_1}}\Bigl[ -c(x,\xi)\ind_{\{\tau_0=0\}}\ind_{\{\rho_0=+\infty\}}+\chi(x,\eta)\ind_{\{\rho_0=0\}}\\
&+V^-\bigl(x+g_\xi(x,\xi)\ind_{\{\tau_0=0\}}\ind_{\{\rho_0=+\infty\}}+g_\eta(x,\eta)\ind_{\{\rho_0=0\}}\bigr)\Bigr],
\end{aligned}
\end{equation*}
therefore
\begin{equation*}
\begin{aligned}
V^-(x)=&\inf_{\underset{\rho_0\in\{0,+\infty\}}{\theta_2(.)\in\Theta_2}}\biggl[\inf_{\eta\in V}\Bigl[\chi(x,\eta)+V^-\bigl(x+g_\eta(x,\eta)\bigr)\Bigr]\ind_{\{\rho_0=0\}}\\
&+\sup_{\underset{\tau_0\in\{0,+\infty\},\;\xi\in U}{\theta_1(.)\in\Theta_1}}\Bigl[-c(x,\xi)\ind_{\{\tau_0=0\}}+V^-\bigl(x+g_\xi(x,\xi)\ind_{\{\tau_0=0\}}\bigr)\Bigr]\ind_{\{\rho_0=+\infty\}}\biggr].
\end{aligned}
\end{equation*}
Since $V^-(x)<\mathcal{H}_{inf}^\chi V^-(x)$, we get
$$V^-(x)=\sup_{\tau_0\in\{0,+\infty\},\;\xi\in U}\Bigl[-c(x,\xi)\ind_{\{\tau_0=0\}}+V^-\bigl(x+g_\xi(x,\xi)\ind_{\{\tau_0=0\}}\bigr)\Bigr].$$
Therefore
$$V^-(x)\geq\sup_{\xi\in U}\Bigl[V^-\bigl(x+g_\xi(x,\xi)\bigr)-c(x,\xi)\Bigr],$$
which completes the proof.
\end{proof}
\begin{lemma}\label{lemma_2}
Assume $\textbf{H}_b$, $\textbf{H}_g$, $\textbf{H}_f$ and $\textbf{H}_{c,\chi}$. Let $x\in\mathbb{R}^n$ and $\phi\in C^1(\mathbb{R}^n)$ be such that
$$\lambda\phi(x)+H^-\bigl(x,D\phi(x)\bigr)=\gamma>0,$$
then there exists $\beta^\gamma\in\mathcal{B}$ a non-anticipative strategy for $player-\eta$ such that, for all $\bigl(\theta_1(.),u\bigr)\in\Theta_1\times\mathcal{U}$ and $t>0$ small enough, we have
$$\int_{0}^{t}\Bigl\{-\lambda \phi\bigl(y_x(s)\bigr)+D\phi\bigl(y_x(s)\bigr).b\bigl(y_x(s);\theta_1(s),\theta_2(s)\bigr)+f\bigl(y_x(s);\theta_1(s),\theta_2(s)\bigr)\Bigr\}\exp(-\lambda s)ds\leq-\frac{\gamma t}{4},$$
where, for $v\in\mathcal{V}$, $y_x(s):=y_x\Bigl(s;\theta_1(.),u,\beta^\gamma\bigl(\theta_1(.),u\bigr)\Bigr)$ and $\beta^\gamma\bigl(\theta_1(.),u\bigr):=\bigl(\theta_2(.),v\bigr)$.\\
A similar result can be obtained for the upper Hamiltonian $H^+$.
\end{lemma}
\begin{proof}
We give only the proof for $H^-$, similarly for $H^+$. Following the results in chapter VIII of reference \cite{BC97}, let $x\in\mathbb{R}^n$, $t>0$ and $\phi\in C^1(\mathbb{R}^n)$ be such that
\begin{equation}\label{equation_13}
\lambda\phi(x)+H^-\bigl(x,D\phi(x)\bigr)=\gamma>0,
\end{equation}
and define for $z\in\mathbb{R}^n$, $\bigl(\theta_1(.),\theta_2(.)\bigr)\in\Theta_1\times\Theta_2$,
$$\Gamma\bigl(z;\theta_1(t),\theta_2(t)\bigr)=\lambda\phi(z)-D\phi(z).b\bigl(z;\theta_1(t),\theta_2(t)\bigr)-f\bigl(z;\theta_1(t),\theta_2(t)\bigr).$$
By (\ref{equation_13}) and the definition of lower Hamiltonian $H^-$ we get $$\inf_{\theta_1\in A}\sup_{\theta_2\in B}\Gamma(x;\theta_1,\theta_2)=\gamma,$$
so for any $\theta_1\in A$ there exists $\theta_2=\theta_2(\theta_1)\in B$ such that $\Gamma\bigl(x;\theta_1,\theta_2\bigr)\geq\gamma$. Since $\theta\rightarrow\Gamma\bigl(x;\theta,\theta_2(t)\bigr)$ is uniformly continuous in $A$, we have in fact
$$\Gamma(x;\zeta,\theta_2)\geq\frac{3\gamma}{4}\;\text{for all}\;\zeta\in B_{r}(\theta_1)\cap A,$$
for some $r=r(\theta_1)$ where $B_{r}\bigl(\theta_1\bigr)$ denotes the open ball of radius $r>0$ centered at $\theta_1$. Since $A$ is a compact subset of $\mathbb{R}^l$, there exist finitely many points $\bigl(\theta_1^1,\theta_1^2,...,\theta_1^n\bigr)$ and $\bigl(r_1,r_2,...,r_n\bigr)$ such that, $\theta_1^i\in A$ and $r_i>0$ for $i=1,2,...,n$ and
$$A\subseteq\cup_{i=1}^{n}B_{r_i}(\theta_1^i),$$
where $r_i:=r_i(\theta_1^i)>0$, and for $\theta_2^i:=\theta_2(\theta_1^i)\in B$
$$\Gamma(x;\zeta,\theta_2^i)\geq\frac{3\gamma}{4}\;\text{for all}\;\zeta\in B_{r_i}(\theta_1^i)\cap A.$$
Next, we define $\psi:\;A\rightarrow B$ by
$$\psi(\theta_1)=\theta_2^k\;\text{if}\;\theta_1\in B_{r_k}(\theta_1^k)\backslash \cup_{i=1}^{k-1}B_{r_i}(\theta_1^i).$$
It is easy to prove that, for any $\theta_1(t)\in A$, $\psi\bigl(\theta_1(t)\bigr)$ is measurable, so we can define $\beta^\gamma\in\mathcal{B}$, a non-anticipative strategy for $player-\eta$, by
$$\beta^\gamma\bigl(\theta_1(t),.\bigr):=\Bigl(\psi\bigl(\theta_1(t)\bigr),.\Bigr).$$
By definition of $\psi$
$$\Gamma\Bigl(x;\theta_1(t),\psi\bigl(\theta_1(t)\bigr)\Bigr)\geq\frac{3\gamma}{4}\;\text{for all}\;\theta_1(t)\in A,$$
and by the continuity of $\Gamma$ and the usual estimate on the trajectories $\bigl\|y_x(t)-y_x(0^-)\bigr\|\leq Ct$ for all $t$ depends on a given constant $C$, its proof is classic (see \cite{BC97}), we deduce that there exists $t>0$ such that
$$\Gamma\Bigl(y_x(s);\theta_1(s),\psi\bigl(\theta_1(s)\bigr)\Bigr)\geq\frac{\gamma}{2}\;\text{for all}\;0\leq s\leq t\;\text{and all}\;\theta_1(s)\in A.$$
Finally we multiply both sides of the last inequality by $\exp(-\lambda s)$ and integrate from $0$ to $t$ to obtain the result for $t$ small enough.
\end{proof}

\subsection{Viscosity Characterization}

Now, we can prove the existence result for the game studied in this paper, i.e., we give the viscosity characterization of the lower and the upper HJBI QVIs (\ref{LowerQVI}) and (\ref{UpperQVI}):
\begin{theorem}\label{ExistenceTh}
Assume $\textbf{H}_b$, $\textbf{H}_g$, $\textbf{H}_f$ and $\textbf{H}_{c,\chi}$. Then the lower value function (\ref{LowerValue}) and the upper value function (\ref{UpperValue}) are viscosity solutions of the lower HJBI QVI (\ref{LowerQVI}) and the upper HJBI QVI (\ref{UpperQVI}), respectively.
\end{theorem}
\begin{proof}
We give only the proof for the lower HJBI QVI (\ref*{LowerQVI}), similarly for the upper HJBI QVI (\ref*{UpperQVI}). A part of this proof is inspired from the results in chapter VIII of reference \cite{BC97}. We first prove the sub-solution property. Let $\phi$ be a function in $C^1(\mathbb{R}^n)$ and $\overline{x}\in\mathbb{R}^n$ be such that $V^{-}-\phi$ achieves a local maximum at $\overline x$ and $V^-(\overline x)=\phi(\overline x)$. If $V^-(\overline x)-\mathcal H_{sup}^c V^-(\overline x)\leq 0$ there is nothing to prove. Otherwise, for $\varepsilon>0$ we assume that $V^-(\overline x)-\mathcal H_{sup}^c V^-(\overline x)\geq\varepsilon>0$, then we proceed by contradiction. Since, from the result of Lemma \ref{lemma_1}, we have $V^-(\overline x)-\mathcal H_{inf}^\chi V^-(\overline x)\leq 0$ we only need to assume that
$$\lambda\phi(\overline x)+H^-\bigl(\overline x,D\phi(\overline x)\bigr)=\gamma>0,$$
then, by the result giving in Lemma \ref{lemma_2}, one can find $\beta^\gamma\in\mathcal{B}$ a non-anticipative strategy for $player-\eta$ such that for all $\theta_1(.)\in\Theta_1$ and $t>0$ small enough
$$\int_{0}^{t}\Bigl\{-\lambda \phi\bigl(y_{\overline x}(s)\bigr)+D\phi\bigl(y_{\overline x}(s)\bigr).b\bigl(y_{\overline x}(s);\theta_1(s),\theta_2(s)\bigr)+f\bigl(y_{\overline x}(s);\theta_1(s),\theta_2(s)\bigr)\Bigr\}\exp(-\lambda s)ds\leq-\frac{\gamma t}{4},$$
where, for $v\in\mathcal{V}$ and any $u\in\mathcal{U}$, $y_{\overline x}(s):=y_{\overline x}\Bigl(s;\theta_1(.),u,\beta^\gamma\bigl(\theta_1(.),u\bigr)\Bigr)$ and $$\beta^\gamma\bigl(\theta_1(.),u\bigr):=\bigl(\theta_2(.),v\bigr),$$
thus,
\begin{equation}\label{equation_14}
\int_{0}^{t}f\bigl(y_{\overline x}(s);\theta_1(s),\theta_2(s)\bigr)\exp(-\lambda s)ds+\exp(-\lambda t)\phi\bigl(y_{\overline x}(t)\bigr)-\phi(\overline x)\leq-\frac{\gamma t}{4}.
\end{equation}
Since $V^{-}-\phi$ has a local maximum at $\overline x$ and $V^{-}(\overline x)=\phi(\overline x)$ we have, for $t$ small enough,
$$\bigl\|y_{\overline x}(t)-\overline x\bigr\|\rightarrow 0,$$
which gives
$$\exp(-\lambda t)\phi\bigl(y_{\overline x}(t)\bigr)-\phi(\overline x)\geq \exp(-\lambda t)V^{-}\bigl(y_{\overline x}(t)\bigr)-V^{-}(\overline x).$$
By plugging this into (\ref{equation_14}) we obtain, for $t=T$ small enough,
\begin{equation*}
\begin{aligned}
\inf_{\beta\in\mathcal{B}}\sup_{\bigl(\theta_1(.),u\bigr)\in\Theta_1\times\mathcal{U}}\biggl\{&\int_0^T f\bigl(y_{\overline x}(t);\theta_1(t),\theta_2(t)\bigr)\exp(-\lambda t)dt\\
&+V^-\biggl(y_{\overline x}\Bigl(T;\theta_1(.),u,\beta\bigl(\theta_1(.),u\bigr)\Bigr)\biggr)\exp(-\lambda T)\biggr\}-V^{-}(\overline x)\leq-\frac{\gamma T}{4}<0,
\end{aligned}
\end{equation*}
where $\beta\bigl(\theta_1(.),u\bigr):=\bigl(\theta_2(.),v\bigr)$. Which, without loss of generality when $T<\tau_0\wedge\rho_0$, is a contradiction to the DPP (\ref{equation_2}), then $V^{-}$ is a viscosity sub-solution of the lower HJBI QVI (\ref*{LowerQVI}).
\par Next, we show the super-solution property. Let $\phi$ be a function in $C^1(\mathbb{R}^n)$ and $\underline{x}\in\mathbb{R}^n$ be such that $V^{-}-\phi$ achieves a local minimum in $B_{\delta}(\underline x)$, where $B_{\delta}(\underline x)$ is the open ball of center $\underline x$ and radius $\delta>0$, and $V^-(\underline x)=\phi(\underline x)$. Now we suppose, for $\varepsilon>0$, that $V^-(\underline x)-\mathcal H_{inf}^\chi V^-(\underline x)<\varepsilon<0$. Then, without loss of generality, we can assume that $V^-(\underline x)-\mathcal H_{inf}^\chi V^-(\underline x)<\varepsilon<0$ on $B_\delta(\underline x)$, then from the result of Lemma \ref{lemma_1} we deduce $V^-(\underline x)-\mathcal H_{sup}^c V^-(\underline x)\geq 0$.
Next, we define
$$t'=\inf\bigl\{t\geq 0: y_{\underline{x}}(t)\notin B_\delta(\underline x)\bigr\}.$$
We let $0<t\leq t'$ and we proceed by contradiction. Assuming that
$$\lambda\phi(\underline x)+H^-\bigl(\underline x,D\phi(\underline x)\bigr)=-\gamma<0,$$
then, by the definition of $H^{-}$, one can find $\alpha^\gamma\in\mathcal{A}$ a non-anticipative strategy for $player-\xi$ such that, for all $\theta_2(.)\in\Theta_2$,
$$\lambda \phi(\underline x)-D\phi(\underline x).b\bigl(\underline x;\theta_1(s),\theta_2(s)\bigr)-f\bigl(\underline x;\theta_1(s),\theta_2(s)\bigr)\leq-\gamma,$$
where, for $u\in\mathcal{U}$ and any $v\in\mathcal{V}$,
$$\alpha^\gamma\bigl(\theta_2(.),v\bigr):=\bigl(\theta_1(.),u\bigr),$$
thus, for $t$ small enough and any $\beta\in\mathcal{B}$
$$\lambda\phi\bigl(y_{\underline x}(s)\bigr)-D\phi\bigl(y_{\underline x}(s)\bigr).b\bigl(y_{\underline x}(s);\theta_1(s),\theta_2(s)\bigr)-f\bigl(y_{\underline x}(s);\theta_1(s),\theta_2(s)\bigr)\leq-\frac{\gamma}{2},$$
where, $0\leq s\leq t$ and for $v\in\mathcal{V}$ and any $u\in\mathcal{U}$, $y_{\underline x}(s):=y_{\underline x}\Bigl(s;\theta_1(.),u,\beta\bigl(\theta_1(.),u\bigr)\Bigr)$ and $$\beta\bigl(\theta_1(.),u\bigr):=\bigl(\theta_2(.),v\bigr).$$
Now we multiply both sides of the last inequality by $\exp(-\lambda s)$ and integrate from $0$ to $t$ to obtain
\begin{equation}\label{equation_15}
\phi(\underline x)-\exp(-\lambda t)\phi\bigl(y_{\underline x}(t)\bigr)-\int_{0}^{t}f\bigl(y_{\underline x};\theta_1(s),\theta_2(s)\bigr)\exp(-\lambda s)ds\leq-\frac{\gamma t}{4}.
\end{equation}
Since $V^{-}-\phi$ has a local minimum at $\underline x$ and $V^{-}(\underline x)=\phi(\underline x)$ we have for $t$ small enough
$$\bigl\|y_{\underline x}(t)-\underline x\bigr\|\rightarrow 0,$$
which gives
$$\exp(-\lambda t)\phi\bigl(y_{\underline x}(t)\bigr)-\phi(\underline x)\leq\exp(-\lambda t)V^{-}\bigl(y_{\underline x}(t)\bigr)-V^{-}(\underline x),$$
thus
$$\exp(-\lambda t)V^{-}\bigl(y_{\underline x}(t)\bigr)+\int_{0}^{t}f\bigl(y_{\underline x};\theta_1(s),\theta_2(s)\bigr)\exp(-\lambda s)ds\geq\frac{\gamma t}{2}+V^{-}(\underline x).$$
By plugging this into (\ref{equation_15}), for $t=T$ small enough, we obtain
\begin{equation*}
\begin{aligned}
\inf_{\beta\in\mathcal{B}}\sup_{\bigl(\theta_1(.),u\bigr)\in\Theta_1\times\mathcal{U}}\biggl\{&\int_0^T f\bigl(y_{\underline x}(t);\theta_1(t),\theta_2(t)\bigr)\exp(-\lambda t)dt\\
&+V^-\biggl(y_{\underline x}\Bigl(T;\theta_1(.),u,\beta\bigl(\theta_1(.),u\bigr)\Bigr)\biggr)\exp(-\lambda T)\biggr\}-V^{-}(\underline x)>0,
\end{aligned}
\end{equation*}
where $\beta\bigl(\theta_1(.),u\bigr):=\bigl(\theta_2(.)\bigr)$. Which, without loss of generality when $T<\tau_0\wedge\rho_0$, is a contradiction to the DPP (\ref{equation_2}), then $V^{-}$ is a viscosity super-solution of the lower HJBI QVI (\ref*{LowerQVI}). The proof is now complete.
\end{proof}
Analogously we introduce the two following HJBI QVIs related, respectively, to the lower Hamiltonian $H^-$ and the upper Hamiltonian $H^+$:
\begin{equation}\label{lowerQVImax}\tag{L$_{\max}$}
\max\biggl\{\min\Bigl[\lambda v(x)+H^-\bigl(x,Dv(x)\bigr),v(x)-\mathcal{H}_{sup}^c v(x)\Bigr],v(x)-\mathcal{H}_{inf}^\chi v(x)\biggr\}=0,
\end{equation}
and
\begin{equation}\label{upperQVImin}\tag{U$_{\min}$}
\min\biggl\{\max\Bigl[\lambda v(x)+H^+\bigl(x,Dv(x)\bigr),v(x)-\mathcal{H}_{inf}^\chi v(x)\Bigr],v(x)-\mathcal{H}_{sup}^c v(x)\biggr\}=0.
\end{equation}
Similarly, we prove the following theorem:
\begin{theorem}\label{Theorem4}
Assume $\textbf{H}_b$, $\textbf{H}_g$, $\textbf{H}_f$ and $\textbf{H}_{c,\chi}$. Then the lower value function (\ref{LowerValue}) and the upper value function (\ref{UpperValue}) are viscosity solutions of the HJBI QVI (\ref{lowerQVImax}) and the HJBI QVI (\ref{upperQVImin}), respectively.\qed
\end{theorem}

\section{Uniqueness of the Viscosity Solution}\label{Sect.4}

In this section we prove the main result of the paper, namely uniqueness wich is inspired from \cite{ElBB10, El17}. First, in Proposition \ref{proposition3}, we give a new formulation of the definition of viscosity solution for the HJBI QVIs (\ref{LowerQVI}) and (\ref{UpperQVI}), which, combined with Lemma \ref{lemma3} below, will be useful to prove the comparison result of Theorem \ref{compTh} hereafter. Next, we conclude in Corollary \ref{corollary2}.
\begin{proposition}{(Viscosity Solution)}\label{proposition3}
A continuous function $V$ in $\mathbb{R}^n$ is a viscosity solution of the lower HJBI QVI (\ref{LowerQVI}) (resp. upper HJBI QVI (\ref{UpperQVI})) if and only if the two following properties hold:
\begin{enumerate}[i.]
\item\underline{Viscosity sub-solution}: For any $\overline x\in\mathbb{R}^n$ and any function $\phi\in C^{1}(\mathbb{R}^n)$ such that $\overline x$ is a local maximum point of $V-\phi$, we have:
\begin{equation*}
\begin{aligned}
\lambda V(\overline x)\leq \max_{i\in\{0,1\}}\biggl\{(1-i)\min_{j\in\{0,1\}}\Bigl[&(1-j)\sup_{\theta_1\in A}\inf_{\theta_2\in B}\bigl\{D\phi(\overline x).b(\overline x;\theta_1,\theta_2)+f(\overline x;\theta_1,\theta_2)\bigr\}\\
&+j\lambda\mathcal{H}_{inf}^\chi V(\overline x)\Bigr]+i\lambda\mathcal{H}_{sup}^c V(\overline x)\biggr\}
\end{aligned}
\end{equation*}
\begin{equation*}
\begin{aligned}
\biggl(\text{resp}.\;\lambda V(\overline x)\leq \min_{i\in\{0,1\}}\biggl\{(1-i)\max_{j\in\{0,1\}}\Bigl[&(1-j)\inf_{\theta_2\in B}\sup_{\theta_1\in A}\bigl\{D\phi(\overline x).b(\overline x;\theta_1,\theta_2)+f(\overline x;\theta_1,\theta_2)\bigr\}\\
&+j\lambda\mathcal{H}_{sup}^c V(\overline x)\Bigr]+i\lambda\mathcal{H}_{inf}^\chi V(\overline x)\biggr\}\biggr);
\end{aligned}
\end{equation*}
\item\underline{Viscosity super-solution}: For any $\underline x\in\mathbb{R}^n$ and any function $\phi\in C^{1}(\mathbb{R}^n)$ such that $\underline x$ is a local minimum point of $V-\phi$, we have:
\begin{equation*}
\begin{aligned}
\lambda V(\underline x)\geq \max_{i\in\{0,1\}}\biggl\{(1-i)\min_{j\in\{0,1\}}\Bigl[&(1-j)\sup_{\theta_1\in A}\inf_{\theta_2\in B}\bigl\{D\phi(\underline x).b(\underline x;\theta_1,\theta_2)+f(\underline x;\theta_1,\theta_2)\bigr\}\\
&+j\lambda\mathcal{H}_{inf}^\chi V(\underline x)\Bigr]+i\lambda\mathcal{H}_{sup}^c V(\underline x)\biggr\}
\end{aligned}
\end{equation*}
\begin{equation*}
\begin{aligned}
\biggl(\text{resp}.\;\lambda V(\underline x)\geq \min_{i\in\{0,1\}}\biggl\{(1-i)\max_{j\in\{0,1\}}\Bigl[&(1-j)\inf_{\theta_2\in B}\sup_{\theta_1\in A}\bigl\{D\phi(\underline x).b(\underline x;\theta_1,\theta_2)+f(\underline x;\theta_1,\theta_2)\bigr\}\\
&+j\lambda\mathcal{H}_{sup}^c V(\underline x)\Bigr]+i\lambda\mathcal{H}_{inf}^\chi V(\underline x)\biggr\}\biggr).
\end{aligned}
\end{equation*}
\end{enumerate}
\end{proposition}
\begin{proof}
This proof is inspired from \cite{El17}. We give only the proof for the lower HJBI QVI (\ref*{LowerQVI}), similarly for the upper HJBI QVI (\ref*{UpperQVI}). For any positive numbers $a,b,a^\prime$ and $b^\prime$, solving a QVI of the form $$\min\Bigl\{\max\bigl[A,B\bigr],C\Bigr\}=0$$
is equivalent to solve the following equation
\begin{equation}\label{equivQVI}
\min_{i\in\{0,1\}}\Bigl\{(1-i)a\max_{j\in\{0,1\}}\bigl[(1-j)a^\prime A+jb^\prime B\bigr]+ibC\Bigr\}=0.
\end{equation}
The same for the inequalities
$$\min\Bigl\{\max\bigl[A,B\bigr],C\Bigr\}\leq 0,\;\text{and}\;\min\Bigl\{\max\bigl[A,B\bigr],C\Bigr\}\geq 0.$$
We use (\ref{equivQVI}), for $a=a^\prime=1$ and $b=b^\prime=\lambda$, to rewrite the lower HJBI QVI (\ref*{LowerQVI}) as follows
\begin{equation*}
\begin{aligned}
\min_{i\in\{0,1\}}\biggl\{(1-i)\max_{j\in\{0,1\}}\Bigl[&(1-j)\inf_{\theta_1\in A}\sup_{\theta_2\in B}\bigl\{\lambda v(x)-Dv(x).b(x;\theta_1,\theta_2)-f(x;\theta_1,\theta_2)\bigr\}\\
&+j\lambda\bigl(v(x)-\mathcal{H}_{inf}^\chi v(x)\bigr)\Bigr]+i\lambda\bigl(v(x)-\mathcal{H}_{sup}^c v(x)\bigr)\biggr\}=0,
\end{aligned}
\end{equation*}
where $v$ being a continuous function in $\mathbb R^n$ and $x$ an element of $\mathbb R^n$. We then get
\begin{equation*}
\begin{aligned}
\min_{i\in\{0,1\}}\biggl\{(1-i)\max_{j\in\{0,1\}}\Bigl[&\lambda v(x)-j\lambda v(x)+(1-j)\inf_{\theta_1\in A}\sup_{\theta_2\in B}\bigl\{-Dv(x).b(x;\theta_1,\theta_2)-f(x;\theta_1,\theta_2)\bigr\}\\
&+j\lambda\bigl(v(x)-\mathcal{H}_{inf}^\chi v(x)\bigr)\Bigr]+i\lambda\bigl(v(x)-\mathcal{H}_{sup}^c v(x)\bigr)\biggr\}=0,
\end{aligned}
\end{equation*}
from which it follows that
\begin{equation*}
\begin{aligned}
\min_{i\in\{0,1\}}\biggl\{(1-i)\max_{j\in\{0,1\}}\Bigl[&\lambda v(x)-(1-j)\sup_{\theta_1\in A}\inf_{\theta_2\in B}\bigl\{Dv(x).b(x;\theta_1,\theta_2)+f(x;\theta_1,\theta_2)\bigr\}\\
&-j\lambda\mathcal{H}_{inf}^\chi v(x)\Bigr]+i\lambda\bigl(v(x)-\mathcal{H}_{sup}^c v(x)\bigr)\biggr\}=0.
\end{aligned}
\end{equation*}
Then we deduce
\begin{equation*}
\begin{aligned}
\max_{i\in\{0,1\}}\biggl\{(1-i)\min_{j\in\{0,1\}}\Bigl[&-\lambda v(x)+(1-j)\sup_{\theta_1\in A}\inf_{\theta_2\in B}\bigl\{Dv(x).b(x;\theta_1,\theta_2)+f(x;\theta_1,\theta_2)\bigr\}\\
&+j\lambda\mathcal{H}_{inf}^\chi v(x)\Bigr]-i\lambda\bigl(v(x)-\mathcal{H}_{sup}^c v(x)\bigr)\biggr\}=0,
\end{aligned}
\end{equation*}
thus
\begin{equation*}
\begin{aligned}
\max_{i\in\{0,1\}}\biggl\{-\lambda v(x)+(1-i)\min_{j\in\{0,1\}}\Bigl[&(1-j)\sup_{\theta_1\in A}\inf_{\theta_2\in B}\bigl\{Dv(x).b(x;\theta_1,\theta_2)+f(x;\theta_1,\theta_2)\bigr\}\\
&+j\lambda\mathcal{H}_{inf}^\chi v(x)\Bigr]+i\lambda\mathcal{H}_{sup}^c v(x)\biggr\}=0.
\end{aligned}
\end{equation*}
Finally we deduce the following expression of the lower HJBI QVI (\ref*{LowerQVI})
\begin{equation}\label{equation_17}
\begin{aligned}
\lambda v(x)=\max_{i\in\{0,1\}}\biggl\{(1-i)\min_{j\in\{0,1\}}\Bigl[&(1-j)\sup_{\theta_1\in A}\inf_{\theta_2\in B}\bigl\{Dv(x).b(x;\theta_1,\theta_2)+f(x;\theta_1,\theta_2)\bigr\}\\
&+j\lambda\mathcal{H}_{inf}^\chi v(x)\Bigr]+i\lambda\mathcal{H}_{sup}^c v(x)\biggr\},
\end{aligned}
\end{equation}
which, using the definition of the viscosity solution for (\ref{equation_17}), completes the proof.
\end{proof}
Next, we give the following useful lemma for which the proof is obvious:
\begin{lemma}\label{lemma3}
If a continuous function $v$ is a viscosity solution to the lower HJBI QVI (\ref{LowerQVI}) (resp. upper HJBI QVI (\ref{UpperQVI})) then for any $0<\mu<1$, $\mu v$ is a viscosity solution to the following QVI:
$$\min\biggl\{\max\Bigl[\lambda v(x)+H_\mu^-\bigl(x,Dv(x)\bigr),v(x)-\mathcal{H}_{inf}^{\chi,\mu} v(x)\Bigr],v(x)-\mathcal{H}_{sup}^{c,\mu} v(x)\biggr\}=0$$
$$\biggl(\text{resp}.\;\max\biggl\{\min\Bigl[\lambda v(x)+H_\mu^+\bigl(x,Dv(x)\bigr),v(x)-\mathcal{H}_{sup}^{c,\mu} v(x)\Bigr],v(x)-\mathcal{H}_{inf}^{\chi,\mu} v(x)\biggr\}=0\biggr),$$
where
$$\mathcal{H}_{inf}^{\chi,\mu} v(x):=\inf_{\eta\in V}\Bigl[v\bigl(x+g_\eta(x,\eta)\bigr)+\mu\chi(x,\eta)\Bigr],$$
$$\mathcal{H}_{sup}^{c,\mu} v(x):=\sup_{\xi\in U}\Bigl[v\bigl(x+g_\xi(x,\xi)\bigr)-\mu c(x,\xi)\Bigr]$$
and
$$H_\mu^-\bigl(x,Dv(x)\bigr):=\inf_{\theta_1\in A}\sup_{\theta_2\in B}\bigl(-Dv(x).b(x;\theta_1,\theta_2)-\mu f(x;\theta_1,\theta_2)\bigr)$$
$$\Bigl(\text{resp}.\;H_\mu^+\bigl(x,Dv(x)\bigr):=\sup_{\theta_2\in B}\inf_{\theta_1\in A}\bigl(-Dv(x).b(x;\theta_1,\theta_2)-\mu f(x;\theta_1,\theta_2)\bigr)\Bigr).$$ \qed
\end{lemma}
Now we are ready to establish the following comparison theorem which is inspired from \cite{ElBB10, El17} and leads us to the uniqueness result for the HJBI QVIs (\ref{LowerQVI}) and (\ref{UpperQVI}):
\begin{theorem}{(Comparison Theorem)}\label{compTh}
Assume $\textbf{H}_b$, $\textbf{H}_g$, $\textbf{H}_f$ and $\textbf{H}_{c,\chi}$. If $u$ is a bounded and uniformly continuous viscosity sub-solution of the lower HJBI QVI (\ref{LowerQVI}) and $v$ is a bounded and uniformly continuous viscosity super-solution of the lower HJBI QVI (\ref{LowerQVI}), then for all $x$ in $\mathbb{R}^n$ we have $u(x)\leq v(x)$. The same result holds true for the upper HJBI QVI (\ref{UpperQVI}).
\end{theorem}
\begin{proof}
We give only the proof for the lower HJBI QVI (\ref*{LowerQVI}), similarly for the upper HJBI QVI (\ref*{UpperQVI}). Let $u$ and $v$ be, respectively, a bounded and uniformly continuous viscosity sub-solution and super-solution to the lower HJBI QVI (\ref*{LowerQVI}). For all $0<\mu<1$, applying Proposition \ref{proposition3} and Lemma \ref{lemma3}, we get that $\mu u$ is a viscosity sub-solution to the following QVI:
\begin{equation}\label{equation_18}
\begin{aligned}
\lambda V(x)=\max_{i\in\{0,1\}}\biggl\{(1-i)\min_{j\in\{0,1\}}\Bigl[&(1-j)\sup_{\theta_1\in A}\inf_{\theta_2\in B}\bigl\{DV(x).b(x;\theta_1,\theta_2)+\mu f(x;\theta_1,\theta_2)\bigr\}\\
&+j\lambda\mathcal{H}_{inf}^{\chi,\mu} V(x)\Bigr]+i\lambda\mathcal{H}_{sup}^{c,\mu} V(x)\biggr\},
\end{aligned}
\end{equation}
where $\mathcal{H}_{inf}^{\chi,\mu}$ and $\mathcal{H}_{sup}^{c,\mu}$ are defined as in Lemma \ref{lemma3} and $V:\mathbb{R}^n\rightarrow\mathbb{R}$. Let us assume that $$M=\sup_{x\in\mathbb{R}^n}\bigl(u(x)-v(x)\bigr)>0,$$ if it is not the case, i.e., $M\leq 0$, then the proof is finished. Then, if $\|u\|_{\infty}\neq 0$ we let $1-M/2\|u\|_{\infty}\leq\mu<1$ to get $$M_\mu=\sup_{x\in\mathbb{R}^n}\bigl(\mu u(x)-v(x)\bigr)>0,$$
otherwise, the fact that $M_\mu>0$ is obvious. The proof will now be divided into three steps:\\
\par \textbf{Step 1.} Let $\varepsilon>0$, $\beta>0$ and consider for any $x,y\in\mathbb{R}^n$ the following test function:
$$\psi_{\mu,\varepsilon,\beta}(x,y)=\mu u(x)-v(y)-\frac{\|x-y\|^2}{\varepsilon^2}-\beta\bigl(\|x\|^2+\|y\|^2\bigr).$$
Let $(x_m,y_m)$ be a maximum point of $\psi_{\mu,\varepsilon,\beta}$ which exists, since this is a continuous function going to infinity when $x$ or $y$ does, and denote $$M_{\psi_{\mu,\varepsilon,\beta}}=\psi_{\mu,\varepsilon,\beta}(x_m,y_m).$$
By definition of $(x_m,y_m)$ we have for all $x,y\in\mathbb{R}^n$,
\begin{equation}\label{equation_19}
\mu u(x_m)-v(y_m)-\frac{\|x_m-y_m\|^2}{\varepsilon^2}-\beta\bigl(\|x_m\|^2+\|y_m\|^2\bigr)\geq \mu u(x)-v(y)-\frac{\|x-y\|^2}{\varepsilon^2}-\beta\bigl(\|x\|^2+\|y\|^2\bigr).
\end{equation}
\begin{itemize}
\item Firstly, we use (\ref{equation_19}) with $y=y_m$ to get that $x_m$ is a maximal point of $\mu u(x)-\phi_u(x)$, where $$\phi_u(x)=\frac{\|x-y_m\|^2}{\varepsilon^2}+\beta\|x\|^2,$$
then, since $\mu u$ is viscosity sub-solution of (\ref{equation_18}), we get
\begin{equation}\label{equation_20}
\begin{aligned}
\lambda\mu u(x_m)\leq\max_{i\in\{0,1\}}\biggl\{(1-i)\min_{j\in\{0,1\}}\biggl[&(1-j)\sup_{\theta_1\in A}\inf_{\theta_2\in B}\Bigl\{\Bigl\langle\frac{2\|x_m-y_m\|}{\varepsilon^2}+2\beta x_m,b(x_m;\theta_1,\theta_2)\Bigr\rangle\\
&+\mu f(x_m;\theta_1,\theta_2)\Bigr\}+j\lambda\mathcal{H}_{inf}^{\chi,\mu}\mu u(x_m)\biggr]+i\lambda\mathcal{H}_{sup}^{c,\mu}\mu u(x_m)\biggr\}.
\end{aligned}
\end{equation}
\item Secondly, we use (\ref{equation_19}) with $x=x_m$ to get that $y_m$ is a minimal point of $v(y)-\phi_v(y)$, where $$\phi_v(y)=-\frac{\|x_m-y\|^2}{\varepsilon^2}-\beta\|y\|^2,$$
then, since $v$ is viscosity super-solution of lower HJBI QVI (\ref{LowerQVI}), by applying Proposition \ref{proposition3} we get
\begin{equation}\label{equation_21}
\begin{aligned}
\lambda v(y_m)\geq\max_{i\in\{0,1\}}\biggl\{(1-i)\min_{j\in\{0,1\}}\biggl[&(1-j)\sup_{\theta_1\in A}\inf_{\theta_2\in B}\Bigl\{\Bigl\langle \frac{2\|x_m-y_m\|}{\varepsilon^2}-2\beta y_m,b(y_m;\theta_1,\theta_2)\Bigr\rangle\\
&+f(x_m;\theta_1,\theta_2)\Bigr\}+j\lambda\mathcal{H}_{inf}^{\chi} v(y_m)\biggr]+i\lambda\mathcal{H}_{sup}^{c} v(y_m)\biggr\}.
\end{aligned}
\end{equation}
\end{itemize}
Hence, using inequalities (\ref{equation_20}) and (\ref{equation_21}), we get
\begin{equation*}
\begin{aligned}
\lambda\Bigl(\mu u(x_m)-v(y_m)\Bigr)\leq& \max_{i\in\{0,1\}}\biggl\{(1-i)\min_{j\in\{0,1\}}\biggl[(1-j)\sup_{\theta_1\in A}\inf_{\theta_2\in B}\Bigl\{\Bigl\langle \frac{2\|x_m-y_m\|}{\varepsilon^2}+2\beta x_m,\\
&b(x_m;\theta_1,\theta_2)\Bigr\rangle+\mu f(x_m;\theta_1,\theta_2)\Bigr\}+j\lambda\mathcal{H}_{inf}^{\chi,\mu}\mu u(x_m)\biggr]+i\lambda\mathcal{H}_{sup}^{c,\mu}\mu u(x_m)\biggr\}\\
&+\min_{i\in\{0,1\}}\biggl\{(1-i)\max_{j\in\{0,1\}}\biggl[(1-j)\inf_{\theta_1\in A}\sup_{\theta_2\in B}\Bigl\{-\Bigl\langle \frac{2\|x_m-y_m\|}{\varepsilon^2}-2\beta y_m,\\
&b(y_m;\theta_1,\theta_2)\Bigr\rangle-f(x_m;\theta_1,\theta_2)\Bigr\}-j\lambda\mathcal{H}_{inf}^{\chi} v(y_m)\biggr]-i\lambda\mathcal{H}_{sup}^{c} v(y_m)\biggr\},
\end{aligned}
\end{equation*}
then
\begin{equation*}
\begin{aligned}
\lambda\Bigl(\mu u(x_m)-v(y_m)\Bigr)\leq& \min_{i\in\{0,1\}}\biggl\{(1-i)\max_{j\in\{0,1\}}\biggl[(1-j)\inf_{\theta_1\in A}\sup_{\theta_2\in B}\Bigl\{\Bigl\langle \frac{2\|x_m-y_m\|}{\varepsilon^2},\\
&b(x_m;\theta_1,\theta_2)-b(y_m;\theta_1,\theta_2)\Bigr\rangle+2\beta\bigl\langle x_m,b(x_m;\theta_1,\theta_2)\bigr\rangle+2\beta\bigl\langle y_m,b(y_m;\theta_1,\theta_2)\bigr\rangle\\
&+\mu f(x_m;\theta_1,\theta_2)-f(y_m;\theta_1,\theta_2)\Bigr\}+j\lambda\bigl(\mathcal{H}_{inf}^{\chi,\mu}\mu u(x_m)-\mathcal{H}_{inf}^{\chi} v(y_m)\bigr)\biggr]\\
&+i\lambda\bigl(\mathcal{H}_{sup}^{c,\mu}\mu u(x_m)-\mathcal{H}_{sup}^{c} v(y_m)\bigr)\biggr\}.
\end{aligned}
\end{equation*}
Thus
\begin{equation}\label{equation_22}
\begin{aligned}
\lambda\Bigl(\mu u(x_m)-v(y_m)\Bigr)\leq \min\biggl\{\max\biggl[& 2C_b\frac{\|x_m-y_m\|^2}{\varepsilon^2}+2\beta\|b\|_\infty\Bigl(\|x_m\|+\|y_m\|\Bigr)+(1-\mu)\|f\|_\infty,\\
&\lambda\Bigl(\mathcal{H}_{inf}^{\chi,\mu} \mu u(x_m)-\mathcal{H}_{inf}^{\chi,\mu} \mu u(y_m)+\bigl\|\bigl(\mathcal{H}_{inf}^{\chi,\mu}\mu u-\mathcal{H}_{inf}^{\chi}v\bigr)^+\bigr\|_\infty\Bigr)\biggr],\\
&\lambda\Bigl(\mathcal{H}_{sup}^{c,\mu} \mu u(x_m)-\mathcal{H}_{sup}^{c,\mu} \mu u(y_m)+\bigl\|\bigl(\mathcal{H}_{sup}^{c,\mu}\mu u-\mathcal{H}_{sup}^{c}v\bigr)^+\bigr\|_\infty\Bigr)\biggr\}.
\end{aligned}
\end{equation}
In the last two steps we investigate the equation in the right-hand side of (\ref{equation_22}), step 2 is devoted to the first term of the equation whereas step 3 concerns the obstacles.\\
\par\textbf{Step 2.} We will prove, in the following, that
\begin{equation}\label{equation_23}
\forall\eta>0,\;\exists\varepsilon_0>0,\;\beta_0>0,\;\forall\varepsilon\leq\varepsilon_0,\;\beta\leq\beta_0:\;\frac{\|x_m-y_m\|^2}{\varepsilon^2}+\beta\bigl(\|x_m\|^2+\|y_m\|^2\bigr)\leq\eta.
\end{equation}
We use inequality (\ref{equation_19}) for $x=y$ then
$$M_{\psi_{\mu,\varepsilon,\beta}}\geq \mu u(x)-v(x)-2\beta\|x\|^2,$$
and we let $\sup_{x\in\mathbb{R}^n}\bigl(\mu u(x)-v(x)\bigr)$ be reached, within $\delta>0$ arbitrary small, in a point $x^*$,
$$\mu u(x^*)-v(x^*)\geq M_\mu-\delta.$$
We choose $\delta$ and $\beta$ such that $M_\mu-\delta-2\beta\|x^*\|^2>0$, which is possible since $x^*$ depends only on $\delta$. Thus we get
\begin{equation}\label{equation_24}
\begin{aligned}
M_{\psi_{\mu,\varepsilon,\beta}}&\geq \mu u(x^*)-v(x^*)-2\beta\|x^*\|^2\\
&\geq M_\mu-\delta-2\beta\|x^*\|^2\\
&>0.
\end{aligned}
\end{equation}
Let $r^2=\mu\|u\|_\infty+\|v\|_\infty$, then
$$\|u\|_\infty\leq M_{\psi_{\mu,\varepsilon,\beta}}\leq r^2-\frac{\|x_m-y_m\|^2}{\varepsilon^2}-\beta\bigl(\|x_m\|^2+\|y_m\|^2\bigr),$$
it follows that
\begin{equation}\label{equation_25}
\|x_m-y_m\|\leq r\varepsilon.
\end{equation}
Therefore we introduce the following increasing function:
$$m(w)=\sup_{\|x-y\|\leq w}|v(x)-v(y)|,$$
then, combining with (\ref{equation_25}), we obtain
$$\mu u(x_m)-v(y_m)=\mu u(x_m)-v(x_m)+v(x_m)-v(y_m)\leq M_\mu+m(r\varepsilon).$$
From the definition of $M_{\psi_{\mu,\varepsilon,\beta}}$ and (\ref{equation_24}) we get
$$M_\mu-\delta-2\beta\|x^*\|^2\leq M_{\psi_{\mu,\varepsilon,\beta}}\leq M_\mu+m(r\varepsilon)-\frac{\|x_m-y_m\|^2}{\varepsilon^2}-\beta\bigl(\|x_m\|^2+\|y_m\|^2\bigr),$$
then $$\frac{\|x_m-y_m\|^2}{\varepsilon^2}+\beta\bigl(\|x_m\|^2+\|y_m\|^2\bigr)\leq\delta+2\beta\|x^*\|^2+m(r\varepsilon).$$
Now we choose $\eta<4M_\mu/3$ and we take $\delta=\eta/4$ and $\beta_0=1$ if $\|x^*\|=0$, $\beta_0=\varepsilon/4\|x^*\|^2$ if $\|x^*\|\neq 0$, to get (\ref{equation_23}), the desired inequality. We also get for any $\beta\leq\beta_0$,
\begin{equation}\label{equation_26}
0<M_\mu-\frac{3\eta}{4}\leq M_\mu-\delta-2\beta\|x^*\|^2\leq M_{\psi_{\mu,\varepsilon,\beta}}\leq \mu u(x_m)-v(y_m).
\end{equation}
\par\textbf{Step 3.} We deduce the contradiction. By (\ref{equation_23}), for $\varepsilon\leq\varepsilon_0$ and $\beta\leq\beta_0$ we have
$$2C_b\|x_m-y_m\|^2/\varepsilon^2\leq 2C_b\eta,\;\beta\|x_m\|\leq\sqrt{\beta\eta},\;\text{and}\;\beta\|y_m\|\leq\sqrt{\beta\eta}.$$
Then, for all $\beta\leq\beta_1=\min\bigl\{\beta_0,\eta/\|b\|^2_\infty\bigr\}$, we get $2\beta\|b\|_\infty\bigl(\|x_m\|+\|y_m\|\bigr)\leq 4\eta$. Moreover, for all $\varepsilon\leq\varepsilon_1=\min\bigl\{\varepsilon_0,\sqrt{\eta}/C_f\bigr\}$, we have $C_f\bigl(\|x_m-y_m\|\bigr)\leq\eta$. By Proposition \ref{proposition2}, $\mathcal{H}_{inf}^{\chi,\mu}\mu u$ and $\mathcal{H}_{sup}^{c,\mu}\mu u$ are uniformly continuous, then, tacking into account (\ref{equation_25}), we find $\varepsilon_2\leq\varepsilon_1$ such that for $\varepsilon\leq\varepsilon_2$, $$\mathcal{H}_{inf}^{\chi,\mu} \mu u(x_m)-\mathcal{H}_{inf}^{\chi,\mu} \mu u(y_m)\leq\eta,\;\text{and}\;\mathcal{H}_{sup}^{c,\mu} \mu u(x_m)-\mathcal{H}_{sup}^{c,\mu} \mu u(y_m)\leq\eta.$$
Thus, tacking into account (\ref{equation_22}), we get for all $\varepsilon\leq\varepsilon_2$ and $\beta\leq\beta_1$,
\begin{equation*}
\begin{aligned}
\lambda\Bigl(\mu u(x_m)-v(y_m)\Bigr)\leq \min\biggl\{\max\Bigl[& (1-\mu)\|f\|_\infty,\lambda\bigl\|\bigl(\mathcal{H}_{inf}^{\chi,\mu}\mu u-\mathcal{H}_{inf}^{\chi}v\bigr)^+\bigr\|_\infty\Bigr],\\
&\lambda\bigl\|\bigl(\mathcal{H}_{sup}^{c,\mu}\mu u-\mathcal{H}_{sup}^{c}v\bigr)^+\bigr\|_\infty\biggr\}+(5+2C_b+\lambda)\eta,
\end{aligned}
\end{equation*}
using (\ref{equation_26}) and the fact that $\eta$ is arbitrary we deduce
\begin{equation*}
\begin{aligned}
\lambda\|(\mu u-v)^+\|_\infty\leq\min\biggl\{\max\Bigl[& (1-\mu)\|f\|_\infty,\lambda\bigl\|\bigl(\mathcal{H}_{inf}^{\chi,\mu}\mu u-\mathcal{H}_{inf}^{\chi}v\bigr)^+\bigr\|_\infty\Bigr],\\
&\lambda\bigl\|\bigl(\mathcal{H}_{sup}^{c,\mu}\mu u-\mathcal{H}_{sup}^{c}v\bigr)^+\bigr\|_\infty\biggr\},
\end{aligned}
\end{equation*}
thus
\begin{equation}\label{equation_27}
\lambda\|(\mu u-v)^+\|_\infty\leq\max\Bigl[ (1-\mu)\|f\|_\infty,\lambda\bigl\|\bigl(\mathcal{H}_{inf}^{\chi,\mu}\mu u-\mathcal{H}_{inf}^{\chi}v\bigr)^+\bigr\|_\infty\Bigr].
\end{equation}
Since for all $x\in\mathbb{R}^n$,
\begin{equation}\label{equation_28}
\mathcal{H}_{inf}^{\chi,\mu}\mu u(x)-\mathcal{H}_{inf}^{\chi}v(x)\leq\sup_{\eta\in V}\Bigl(\mu u\bigl(x+g_\eta(x,\eta)\bigr)-v\bigl(x+g_\eta(x,\eta)\bigr)\Bigr)+\sup_{\eta\in V}\bigl((\mu-1)\chi(x,\eta)\bigr).
\end{equation}
We recall that from Assumption $\textbf{H}_{c,\chi}$ for all $x\in\mathbb{R}^n$, $\eta\in V$, $\chi(x,\eta)>0$. Then, since $0<\mu<1$, from (\ref{equation_28}) we get
$$\bigl\|\bigl(\mathcal{H}_{inf}^{\chi,\mu}\mu u-\mathcal{H}_{inf}^{\chi}v\bigr)^+\bigr\|_\infty<\|(\mu u-v)^+\|_\infty.$$
Therefore (\ref{equation_27}) and the last inequality imply
$$\lambda\|(\mu u-v)^+\|_\infty\leq (1-\mu)\|f\|_\infty.$$
Finally, by letting $\mu\rightarrow 1$ and since $f$ is bounded, we obtain $\|(u-v)^+\|_\infty\leq 0$, which leads us to a contradiction and gives the desired comparison, for any $x\in\mathbb{R}^n$, $u(x)\leq v(x)$.
\end{proof}
\begin{corollary}\label{corollary1}
Under Assumptions $\textbf{H}_b$, $\textbf{H}_g$, $\textbf{H}_f$ and $\textbf{H}_{c,\chi}$, the lower HJBI QVI (\ref{LowerQVI}) has a unique bounded and uniformly continuous viscosity solution. The same result holds true for the upper HJBI QVI (\ref{UpperQVI}).
\end{corollary}
\begin{proof}
We give only the proof for the lower HJBI QVI (\ref*{LowerQVI}), similarly for the upper HJBI QVI (\ref*{UpperQVI}). Assume that $v_1$ and $v_2$ are two viscosity solutions to the lower HJBI QVI (\ref*{LowerQVI}). We first use $v_1$ as a bounded and uniformly continuous viscosity sub-solution and $v_2$ as a bounded and uniformly continuous viscosity super-solution and we recall the comparison theorem. Then we change the role of $v_1$ and $v_2$ to
get $v_1(x)=v_2(x)$ for any $x\in\mathbb{R}^n$.
\end{proof}
Next, in the following we give the uniqueness result for the game studied in this paper:
\begin{theorem}
Assume $\textbf{H}_b$, $\textbf{H}_g$, $\textbf{H}_f$, $\textbf{H}_{c,\chi}$ and Isaacs' condition $H^-=H^+$. Both the lower and upper HJBI QVIs (\ref{LowerQVI}) and (\ref{UpperQVI}) admit the value function (\ref{value}) as the unique bounded and uniformly continuous viscosity solution.
\end{theorem}
\begin{proof}
The proof follows immediately from Theorem \ref{ExistenceTh}, Theorem \ref{Theorem4} and Corollary \ref{corollary1} because (\ref*{LowerQVI}) and (\ref*{UpperQVI}) coincide with (\ref*{upperQVImin}) and (\ref*{lowerQVImax}), respectively, if $H^-=H^+$.
\end{proof}
\begin{corollary}\label{corollary2}
Under Assumptions $\textbf{H}_b$, $\textbf{H}_g$, $\textbf{H}_f$, $\textbf{H}_{c,\chi}$ and Isaacs' condition $H^-=H^+$, the lower value function (\ref{LowerValue}) and upper value function (\ref{UpperValue}) coincide and the value function (\ref{value}) $V:=V^-=V^+$ of the infinite horizon, two-player, zero-sum, deterministic, differential game involving continuous and impulse controls is the unique viscosity solution to the lower HJBI QVI (\ref{LowerQVI}) (or, upper HJBI QVI (\ref{UpperQVI})).\qed
\end{corollary}
\section{Conclusion}\label{Sect.5}
We have considered a new class of infinite horizon, two-player, zero-sum, deterministic differential games where each player uses both continuous and impulse controls with discounted payoff. We have studied this class of differential games in viscosity solutions framework under the classical assumptions of Section \ref{Sect.2.2}, the value function (\ref{value}) results to be well-posed. To the best of our knowledge, this is the first characterization of such a class of games, that enjoys a wide range of applications in various fields of engineering, such as mathematical finance (see the aforementioned Example \ref{Example.1}, Section \ref{Sect.2.1}). The class of differential games studied has a form of impulses that depends on nonlinear functions $g_\xi$ and $g_\eta$, and a general costs of impulses, costs that are depending on the system's state $y_x(.)$. In this, our paper differs from and extends many earlier results on zero-sum deterministic impulse controls games.
\par We intend to develop this work in two main directions in the future: First, we aim to give a discrete-time approximation of the Hamilton-Jacobi-Bellman-Isaacs equation (\ref{UpperQVI}) introduced in this work, i.e., prove the existence and uniqueness of a function $v_h$, solution of some approximate equation, which tends to the value function (\ref{value}), when $h$ goes to $0$ (see \cite{El17,S85,S85'}). Second, we wish to characterize a Nash-Equilibrium strategy for the class of differential games we have studied. Hence, by giving some meaningful nonlinear functions ($b^\pi,g_\xi^\pi$ and $g_\eta^\pi$), gain $f^\pi$ and costs ($c^\pi$ and $\chi^\pi$) we might derive a new dynamic portfolio optimization model (see the aforementioned Example \ref{Example.1}, Section \ref{Sect.2.1}).

\end{document}